\documentclass[12pt]{amsart}
\usepackage{amssymb,latexsym, mathtools,tikz}
\usepackage{url}
\usetikzlibrary{arrows,decorations.markings}
\tikzstyle arrowstyle=[scale=1]
\tikzstyle directed=[postaction={decorate,
decoration={markings,mark=at position .65 with {\arrow[arrowstyle]{stealth}}}}]

\usepackage{appendix}
\usepackage{extarrows} 
\usepackage{amscd}
\usepackage{xypic}
\usepackage[utf8]{inputenc}
\usepackage[top=35mm, bottom=25mm, left=30mm, right=30mm]{geometry}
\usepackage{graphicx}

\begin{document}

\title[Shift modules, strongly stable ideals, and their dualities]{Shift modules, strongly stable ideals, and their dualities}
       
\author{Gunnar Fl{\o}ystad}
\address{Matematisk Institutt\\
         Postboks\\
         5020 Bergen}
\email{gunnar@mi.uib.no}


\date{\today}

\begin{abstract}
  We enrich the setting of strongly stable ideals (SSI): We introduce
  shift modules, a module category encompassing SSI's.
  The recently introduced duality on SSI's is
  given an effective conceptual and computational setting.
  We study strongly stable ideals in infinite dimensional polynomial rings,
  where the duality is most natural. Finally a new
  type of resolution for SSI's is introduced. This is the projective resolution
  in the category of shift modules.
\end{abstract}
\maketitle


\theoremstyle{plain}
\newtheorem{theorem}{Theorem}[section]
\newtheorem{corollary}[theorem]{Corollary}
\newtheorem*{main}{Main Theorem}
\newtheorem{lemma}[theorem]{Lemma}
\newtheorem{proposition}[theorem]{Proposition}
\newtheorem{conjecture}[theorem]{Conjecture}

\theoremstyle{definition}
\newtheorem{definition}[theorem]{Definition}
\newtheorem{fact}[theorem]{Fact}
\newtheorem{obs}[theorem]{Observation}
\newtheorem{definisjon}[theorem]{Definisjon}
\newtheorem{problem}[theorem]{Problem}
\newtheorem{condition}[theorem]{Condition}

\theoremstyle{remark}
\newtheorem{notation}[theorem]{Notation}
\newtheorem{remark}[theorem]{Remark}
\newtheorem{example}[theorem]{Example}
\newtheorem{claim}{Claim}
\newtheorem{observation}[theorem]{Observation}


\newcommand{\psp}[1]{{{\bf P}^{#1}}}
\newcommand{\psr}[1]{{\bf P}(#1)}
\newcommand{\op}{{\mathcal O}}
\newcommand{\opw}{\op_{\psr{W}}}

\newcommand{\ini}[1]{\text{in}(#1)}
\newcommand{\gin}[1]{\text{gin}(#1)}
\newcommand{\kr}{{\Bbbk}}
\newcommand{\pd}{\partial}
\newcommand{\vardel}{\partial}
\renewcommand{\tt}{{\bf t}}


\newcommand{\coh}{{{\text{{\rm coh}}}}}


\newcommand{\modv}[1]{{#1}\text{-{mod}}}
\newcommand{\modstab}[1]{{#1}-\underline{\text{mod}}}

\newcommand{\sut}{{}^{\tau}}
\newcommand{\sumit}{{}^{-\tau}}
\newcommand{\til}{\thicksim}

\newcommand{\totp}{\text{Tot}^{\prod}}
\newcommand{\dsum}{\bigoplus}
\newcommand{\dprod}{\prod}
\newcommand{\lsum}{\oplus}
\newcommand{\lprod}{\Pi}

\newcommand{\La}{{\Lambda}}

\newcommand{\sirstj}{\circledast}

\newcommand{\she}{\EuScript{S}\text{h}}
\newcommand{\cm}{\EuScript{CM}}
\newcommand{\cmd}{\EuScript{CM}^\dagger}
\newcommand{\cmri}{\EuScript{CM}^\circ}
\newcommand{\cler}{\EuScript{CL}}
\newcommand{\clerd}{\EuScript{CL}^\dagger}
\newcommand{\clerri}{\EuScript{CL}^\circ}
\newcommand{\gor}{\EuScript{G}}
\newcommand{\cF}{\mathcal{F}}
\newcommand{\cG}{\mathcal{G}}
\newcommand{\cM}{\mathcal{M}}
\newcommand{\cE}{\mathcal{E}}
\newcommand{\cI}{\mathcal{I}}
\newcommand{\cP}{\mathcal{P}}
\newcommand{\cK}{\mathcal{K}}
\newcommand{\cS}{\mathcal{S}}
\newcommand{\cC}{\mathcal{C}}
\newcommand{\cO}{\mathcal{O}}
\newcommand{\cJ}{\mathcal{J}}
\newcommand{\cU}{\mathcal{U}}
\newcommand{\cQ}{\mathcal{Q}}
\newcommand{\cX}{\mathcal{X}}
\newcommand{\cY}{\mathcal{Y}}
\newcommand{\cZ}{\mathcal{Z}}
\newcommand{\cV}{\mathcal{V}}

\newcommand{\mm}{\mathfrak{m}}

\newcommand{\dlim} {\varinjlim}
\newcommand{\ilim} {\varprojlim}

\newcommand{\CM}{\text{CM}}
\newcommand{\Mon}{\text{Mon}}


\newcommand{\Kom}{\text{Kom}}


\newcommand{\EH}{{\mathbf H}}
\newcommand{\res}{\text{res}}
\newcommand{\Hom}{\text{Hom}}
\newcommand{\inhom}{{\underline{\text{Hom}}}}
\newcommand{\Ext}{\text{Ext}}
\newcommand{\Tor}{\text{Tor}}
\newcommand{\ghom}{\mathcal{H}om}
\newcommand{\gext}{\mathcal{E}xt}
\newcommand{\id}{\text{{id}}}
\newcommand{\im}{\text{im}\,}
\newcommand{\codim} {\text{codim}\,}
\newcommand{\resol}{\text{resol}\,}
\newcommand{\rank}{\text{rank}\,}
\newcommand{\lpd}{\text{lpd}\,}
\newcommand{\coker}{\text{coker}\,}
\newcommand{\supp}{\text{supp}\,}
\newcommand{\Ad}{A_\cdot}
\newcommand{\Bd}{B_\cdot}
\newcommand{\Fd}{F_\cdot}
\newcommand{\Gd}{G_\cdot}


\newcommand{\sus}{\subseteq}
\newcommand{\sups}{\supseteq}
\newcommand{\pil}{\rightarrow}
\newcommand{\vpil}{\leftarrow}
\newcommand{\rpil}{\leftarrow}
\newcommand{\lpil}{\longrightarrow}
\newcommand{\inpil}{\hookrightarrow}
\newcommand{\pils}{\twoheadrightarrow}
\newcommand{\projpil}{\dashrightarrow}
\newcommand{\dotpil}{\dashrightarrow}
\newcommand{\adj}[2]{\overset{#1}{\underset{#2}{\rightleftarrows}}}
\newcommand{\mto}[1]{\stackrel{#1}\longrightarrow}
\newcommand{\vmto}[1]{\stackrel{#1}\longleftarrow}
\newcommand{\mtoelm}[1]{\stackrel{#1}\mapsto}
\newcommand{\bihom}[2]{\overset{#1}{\underset{#2}{\rightleftarrows}}}
\newcommand{\eqv}{\Leftrightarrow}
\newcommand{\impl}{\Rightarrow}

\newcommand{\iso}{\cong}
\newcommand{\te}{\otimes}
\newcommand{\into}[1]{\hookrightarrow{#1}}
\newcommand{\ekv}{\Leftrightarrow}
\newcommand{\equi}{\simeq}
\newcommand{\isopil}{\overset{\cong}{\lpil}}
\newcommand{\equipil}{\overset{\equi}{\lpil}}
\newcommand{\ispil}{\isopil}
\newcommand{\vvi}{\langle}
\newcommand{\hvi}{\rangle}
\newcommand{\susneq}{\subsetneq}
\newcommand{\sgn}{\text{sign}}


\newcommand{\xd}{\check{x}}
\newcommand{\ortog}{\bot}
\newcommand{\tL}{\tilde{L}}
\newcommand{\tM}{\tilde{M}}
\newcommand{\tH}{\tilde{H}}
\newcommand{\tvH}{\widetilde{H}}
\newcommand{\tvh}{\widetilde{h}}
\newcommand{\tV}{\tilde{V}}
\newcommand{\tS}{\tilde{S}}
\newcommand{\tT}{\tilde{T}}
\newcommand{\tR}{\tilde{R}}
\newcommand{\tf}{\tilde{f}}
\newcommand{\ts}{\tilde{s}}
\newcommand{\tp}{\tilde{p}}
\newcommand{\tr}{\tilde{r}}
\newcommand{\tfst}{\tilde{f}_*}
\newcommand{\empt}{\emptyset}
\newcommand{\bfa}{{\mathbf a}}
\newcommand{\bfb}{{\mathbf b}}
\newcommand{\bfd}{{\mathbf d}}
\newcommand{\bfl}{{\mathbf \ell}}
\newcommand{\bfx}{{\mathbf x}}
\newcommand{\bfm}{{\mathbf m}}
\newcommand{\bfv}{{\mathbf v}}
\newcommand{\bft}{{\mathbf t}}
\newcommand{\bbfa}{{\mathbf a}^\prime}
\newcommand{\la}{\lambda}
\newcommand{\bfen}{{\mathbf 1}}
\newcommand{\bfe}{{\mathbf e}}
\newcommand{\ep}{\epsilon}
\newcommand{\en}{r}
\newcommand{\tu}{s}
\newcommand{\Sym}{\text{Sym}}

\newcommand{\ome}{\omega_E}

\newcommand{\bevis}{{\bf Proof. }}
\newcommand{\demofin}{\qed \vskip 3.5mm}
\newcommand{\nyp}[1]{\noindent {\bf (#1)}}
\newcommand{\demo}{{\it Proof. }}
\newcommand{\demodone}{\demofin}
\newcommand{\parg}{{\vskip 2mm \addtocounter{theorem}{1}  
                   \noindent {\bf \thetheorem .} \hskip 1.5mm }}

\newcommand{\lcm}{{\text{lcm}}}


\newcommand{\dl}{\Delta}
\newcommand{\cdel}{{C\Delta}}
\newcommand{\cdelp}{{C\Delta^{\prime}}}
\newcommand{\dlst}{\Delta^*}
\newcommand{\Sdl}{{\mathcal S}_{\dl}}
\newcommand{\lk}{\text{lk}}
\newcommand{\lkd}{\lk_\Delta}
\newcommand{\lkp}[2]{\lk_{#1} {#2}}
\newcommand{\del}{\Delta}
\newcommand{\delr}{\Delta_{-R}}
\newcommand{\dd}{{\dim \del}}
\newcommand{\Del}{\Delta}

\renewcommand{\aa}{{\bf a}}
\newcommand{\bb}{{\bf b}}
\newcommand{\cc}{{\bf c}}
\newcommand{\xx}{{\bf x}}
\newcommand{\yy}{{\bf y}}
\newcommand{\zz}{{\bf z}}
\newcommand{\mv}{{\xx^{\aa_v}}}
\newcommand{\mF}{{\xx^{\aa_F}}}

\newcommand{\Symm}{\text{Sym}}
\newcommand{\pnm}{{\bf P}^{n-1}}
\newcommand{\opnm}{{\go_{\pnm}}}
\newcommand{\ompnm}{\omega_{\pnm}}

\newcommand{\pn}{{\bf P}^n}
\newcommand{\hele}{{\mathbb Z}}
\newcommand{\nat}{{\mathbb N}}
\newcommand{\rasj}{{\mathbb Q}}
\newcommand{\bfone}{{\mathbf 1}}

\newcommand{\dt}{\bullet}
\newcommand{\disk}{\scriptscriptstyle{\bullet}}

\newcommand{\cxF}{F_\dt}
\newcommand{\pol}{f}

\newcommand{\Rn}{{\mathbb R}^n}
\newcommand{\An}{{\mathbb A}^n}
\newcommand{\frg}{\mathfrak{g}}
\newcommand{\PW}{{\mathbb P}(W)}

\newcommand{\pos}{{\mathcal Pos}}
\newcommand{\g}{{\gamma}}

\newcommand{\Vaa}{V_0}
\newcommand{\Bp}{B^\prime}
\newcommand{\Bpp}{B^{\prime \prime}}
\newcommand{\bbp}{\mathbf{b}^\prime}
\newcommand{\bbpp}{\mathbf{b}^{\prime \prime}}
\newcommand{\bp}{{b}^\prime}
\newcommand{\bpp}{{b}^{\prime \prime}}

\newcommand{\oLa}{\overline{\Lambda}}
\newcommand{\ov}[1]{\overline{#1}}
\newcommand{\ovv}[1]{\overline{\overline{#1}}}
\newcommand{\tm}{\tilde{m}}
\newcommand{\po}{\bullet}

\newcommand{\surj}[1]{\overset{#1}{\twoheadrightarrow}}
\newcommand{\Supp}{\text{Supp}}

\def\CC{{\mathbb C}}
\def\GG{{\mathbb G}}
\def\ZZ{{\mathbb Z}}
\def\NN{{\mathbb N}}
\def\RR{{\mathbb R}}
\def\OO{{\mathbb O}}
\def\QQ{{\mathbb Q}}
\def\VV{{\mathbb V}}
\def\PP{{\mathbb P}}
\def\EE{{\mathbb E}}
\def\FF{{\mathbb F}}
\def\AA{{\mathbb A}}

\newcommand{\oR}{\overline{R}}
\newcommand{\bfu}{{\mathbf u}}
\newcommand{\nn}{{\mathbf n}}
\newcommand{\oa}{\overline{a}}
\newcommand{\cop}{\text{cop}}
\renewcommand{\op}{\text{op}}
\renewcommand{\mm}{{\mathbf m}}
\newcommand{\ngmi}{\text{neg}}
\newcommand{\up}{\text{up}}
\newcommand{\dw}{\text{down}}
\newcommand{\di}[1]{\hat{#1}}
\newcommand{\diw}[1]{\widehat{#1}}
\newcommand{\bo}{b}
\newcommand{\ub}{u}
\newcommand{\fs}{*}
\newcommand{\ifst}{\infty}
\newcommand{\mon}{{mon}}
\newcommand{\cl}{\text{cl}}
\newcommand{\regr}{\text{reg}}
\newcommand{\reg}{\text{reg}}
\newcommand{\ul}[1]{\underline{#1}}
\renewcommand{\ov}[1]{\overline{#1}}
\newcommand{\bipil}{\leftrightarrow}
\newcommand{\bfc}{{\mathbf c}}
\renewcommand{\mp}{m^\prime}
\newcommand{\np}{n^\prime}
\newcommand{\Mod}{\text{Mod }}
\newcommand{\Sh}{\text{Sh } }
\newcommand{\st}{\text{st}}
\newcommand{\hM}{\tilde{M}}
\newcommand{\hs}{\tilde{s}}
\newcommand{\ee}{\mathbf{e}}
\renewcommand{\dd}{\mathbf{d}}
\renewcommand{\en}{{\mathbf 1}}
\long\def\ignore#1{}
\newcommand{\lex}{{\text{lex}}}
\newcommand{\ordGL}{\succeq_{\lex}}
\newcommand{\ordG}{\succ_{\lex}}
\newcommand{\ordML}{\preceq_{\lex}}
\newcommand{\ordM}{\prec_{\lex}}
\newcommand{\tLa}{\tilde{\Lambda}}
\newcommand{\tGa}{\tilde{\Gamma}}
\newcommand{\STS}{\text{STS}}
\newcommand{\ii}{\mathbf{i}}
\newcommand{\jj}{\mathbf{j}}
\renewcommand{\mod}{\text{ mod}\,}
\newcommand{\shmod}{\texttt{shmod}\,}
\newcommand{\hf}{\underline{f}}
\newcommand{\Glim}{\lim}
\newcommand{\Gcolim}{\colim}
\newcommand{\fm}{f^{\underline{m}}}
\newcommand{\fn}{f^{\underline{n}}}
\newcommand{\gn}{g_{\mathbf{|}n}}
\newcommand{\llin}{\raisebox{1pt}{\scalebox{1}[0.6]{$\mid$}}}
\newcommand{\promap}{\mathrlap{{\hskip 2.8mm}{\llin}}{\lpil}}
\newcommand{\Hilb}{\text{Hilb}}
\newcommand{\Pn}{{\mathbb P}^n}
\newcommand{\Spec}{{\rm{Spec}\,}}
\renewcommand{\ini}{\text{in}}

\ignore{
\makeatletter
\newcommand{\colim@}[2]{%
  \vtop{\m@th\ialign{##\cr
    \hfil$#1\operator@font colim$\hfil\cr
    \noalign{\nointerlineskip\kern1.5\ex@}#2\cr
    \noalign{\nointerlineskip\kern-\ex@}\cr}}%
}
\newcommand{\colim}{%
  \mathop{\mathpalette\colim@{\rightarrowfill@\textstyle}}\nmlimits@
}
\makeatother

\makeatletter
\newcommand{\lim@}[2]{%
  \vtop{\m@th\ialign{##\cr
    \hfil$#1\operator@font lim$\hfil\cr
    \noalign{\nointerlineskip\kern1.5\ex@}#2\cr
    \noalign{\nointerlineskip\kern-\ex@}\cr}}%
}
\renewcommand{\lim}{%
  \mathop{\mathpalette\lim@{\leftarrowfill@\textstyle}}\nmlimits@
}
\makeatother
}
\renewcommand{\lim}{\text{lim}}
\newcommand{\colim}{\text{colim}}
\newcommand{\upa}{\uparrow}

\setcounter{tocdepth}{1} 
\tableofcontents


\section{Introduction}

Strongly stable ideals are somewhat hard to place in the landscape of
mathematics, but let us venture a brief tour.
In algebraic geometry the ideal of every projective variety (char. $0$)
degenerates to such an ideal, \cite{Gal} or see \cite[Sec.15.9]{Ei}.
Also called Borel-fixed ideals,
they are a way to understand and classify components of the Hilbert scheme,
\cite{Re-St, Fl-Ro, Ka-Le},  \cite{Ber-Lel-Rog, Cio-Rog,
  BeCiRo},
\cite{Ram-bor, Sta-hil}.
They are the most degenerate of homogeneous ideals in polynomial
rings $k[x_1, \ldots, x_n]$: The closed orbits for the action of $GL(n)$ on
the Hilbert scheme are precisely the orbits of the strongly stable ideals,
see Appendix \ref{app:closed}. However they being so degenerate, there
is hardly any geometry left in them.

In commutative algebra they have a distinguished resolution, the
Eliahou-Kervaire resolution \cite{EK}.
They occur in the study of Hilbert functions
and Betti numbers of graded ideals \cite{Big-betti}, and in particular in the
proofs of Macaulay's theorems \cite{Gre-ini}.
Algorithms to generate Borel-fixed ideals with given invariants
are given in \cite{Lel-alg} and \cite{Mor-Nag}.

In combinatorics one finds them in shifting theory \cite{Kal-shift},
\cite{Bj-Ka-acta},
minimal growth
of Hilbert functions \cite{He-Hi}, and relations to posets
\cite[Section 6]{FGH}. 
One might even consider them so skeletal and degenerate that they are more
or less numerical objects.
In any case they retain significant
invariants of ideals that degenerate to them, for instance regularity
\cite{Ba-St}.

 They occur in a number of places but always on the
 fringe. Their natural position and effective
 use is however clear. Standard references \cite[Ch.4]{He-Hi} and
 \cite[Ch.2]{Mi-St} have early on specific chapters on them with basic and
significant theory. \cite{Pe} has much the same more distributed,
see also \cite[Ch.15]{Ei}. The most comprehensive treatment may be
\cite{Gre-ini} with many examples and relations to algebraic geometry.

In this article we enrich the setting of strongly stable ideals.
The following new features are studied:
\begin{itemize}
\item[1.] Shift modules: extending strongly stable ideals to a category
  of modules,
\item[2.] Dualities: recently discovered in \cite{Fl-pos} and \cite{Sh-Ya},
  \item[3.] Ambient polynomial ring with infinitely many variables:
    natural setting for the dualities,
  \item[4.] Resolutions: new type of projective resolution
    with new homological invariants.
\end{itemize}

\medskip
\noindent {\bf 1. Shift modules.} Over a polynomial ring
$k[x_1, \ldots, x_n]$ where $k$ is a field, we introduce a category of modules
with shift operations.
This class subsumes strongly stable ideals.

Recall that a monomial ideal $I$ in $k[x_1, \ldots, x_n]$ is
strongly stable if whenever a monomial
$x_j u \in I$ and $i < j$, then $x_i u \in I$.
We may write
$x_iu = s_{i,j}(x_ju)$. So an ideal $I$ is strongly stable when
it is invariant under such shift operations. 

This inspires defining the category of multigraded
shift modules over polynomial rings. Such a module $M$ comes
with shift operations between graded pieces
\[ s_{ij}: M_{\bfd + \ee_j} \to M_{\bfd + \ee_i}. \]
A typical example of a shift module comes from an inclusion of two
strongly stable ideals $I \sus J$: The quotient $J/I$ is a shift module.
More generally quotients of maps between sums of such ideals are shift modules.

We define shift modules in
three steps.
First we define {\it finite} shift modules, Section \ref{sec:finmod}.
Such a module is graded by
a finite set of degrees $\dd \in \NN^{m+1}$ with some fixed total degree.
Then we define shift modules over finite dimensional polynomial rings
$k[x_1, \ldots, x_m]$. We show that finitely generated such modules
come from finite shift modules, Theorem \ref{thm:finpol-exp}.
Lastly we define shift modules over the infinite
dimensional polynomial ring.

\medskip
Given a monomial $m$ let $\langle m \rangle$ to be the smallest
strongly stable ideal containing $m$. For instance $\langle x_1 x_2 x_3
\rangle$ will be the ideal generated by
$x_1^3,x_1^2x_2, x_1x_2^2, x_1^2x_3, x_1x_2x_3$. 

\medskip
\noindent {\bf 2. Dualities.} Recently a duality
on strongly stable ideals was discovered, \cite{Fl-pos} and \cite{Sh-Ya}.
We develop the conceptual framework for this duality, enabling effective
arguments and concrete computations of duals. Moreover we extend this
duality to shift modules.  This is an analog of extending Alexander duality for
squarefree monomial ideals to squarefree modules, \cite{Ya-Sq}.

We get nice formulas such as
the dual of $\langle y_m^n \rangle$ being $\langle x_n^m \rangle$,
 Corollary \ref{cor:st-idealer:xam}, and the following:

\medskip
\noindent {\bf Proposition \ref{pro:settCI}}
{\em The dual of the strongly stable ideal with one generator
$\langle y_{a_1} y_{a_2} \cdots y_{a_n} \rangle$
is the strongly stable ideal with generators
    $\langle x_1^{a_1}, x_2^{a_2}, \ldots, x_n^{a_n} \rangle$.}

\medskip
Furthermore the duality takes sums of ideals to intersections of their
duals, Corollary \ref{cor:st-idealer:sums}.

A complicating aspect with these dualities is that they do not take place
in a fixed finite dimensional polynomial ring, as seen by the above.
Fixing $m$ the dual of $\langle y_m^n \rangle $ is $\langle x_n^m \rangle$,
and the latter requires larger and larger polynomial rings as
$n$ gets larger.
In order to have a full
natural setting we must be in an infinite dimensional polynomial ring. 

\medskip
A basic tool and inspiration for our work here, is that strongly stable ideals
generated in degree $\leq n$ in $k[x_1, \ldots, x_m]$ are in
one-to-one correspondence with poset ideals in $\Hom([m],[n+1])$,
\cite[Section 6]{FGH}. To extend this, 
let $\di{\NN}  = \NN \cup \{ \infty \}$. The poset of order preserving maps
$f : \NN \pil \di{\NN}$, denoted $\Hom(\NN, \di{\NN})$ is a
central object for us, and in \cite[Section 5]{Fl-Pro} we introduced
a {\em topology} on this poset.

\medskip
\noindent {\bf 3. Infinite dimensional polynomial rings.}
Let $k[x_\NN]$ be the infinite dimensional polynomial ring
in the variables $x_1, x_2, \cdots $ indexed by natural numbers.
A basic tool we use is the commutative diagram of bijections:
\begin{equation} \label{eq:intro:LG} \xymatrix{\Hom_S(\NN,\di{\NN}) \ar[dd]^{D} \ar[dr]^{\Lambda} & \\
    & \Mon(x_\NN). \\
    \Hom^L(\NN,\di{\NN}) \ar[ur]^{\Gamma} & 
  }
  \end{equation}
  Here $\Hom_S$ are the bounded maps in $\Hom(\NN, \di{\NN})$, and
  $\Hom^L$ are the maps which eventually take value $\infty$.
  The maps are given by
  \[  \Gamma(f) = \prod_{f(i) < \ifst} x_{f(i)}, \quad
    \Lambda(f) = \prod_{i \geq 1} x_i^{f(i)-f(i-1)}.\]
We show that strongly stable ideals in $k[x_\NN]$ are in
one-to-one correspondence with {\em open poset ideals}
in $\Hom(\NN,\di{\NN})$. 
The duality for strongly stable ideals in $k[x_\NN]$ comes about because
$\Hom(\NN, \di{\NN})$ is a {\em self-dual poset},
Subsection \ref{subsec:posiso} .
Not all strongly stable ideals
in $k[x_\NN]$ have duals. We identify precisely which ideals have.
In particular all {\it finitely generated} strongly stable ideals have duals.

We remark that a recent trend in commutative algebra is to consider
infinite dimensional polynomial rings, \cite{AsHi, NaRo, HiSu, KrLeSn}.
The increasing monoid
$\Hom_{inj}(\NN, \NN)$ of injective order-preserving maps $f : \NN \pil \NN$
is in one-one correspondence with $\Hom(\NN,\NN)$ by mapping
$f$ to $f - \id_\NN + 1$.
The increasing monoid
has been used to study $k[x_\NN]$ in  \cite{NaRo}, \cite{KrLeSn}, see also
\cite{GuSn}.
The use differs however sharply from ours, as $\Hom_{inj}(\NN, \NN)$
is used there to act
on $k[x_\NN]$, while we use $\Hom(\NN, \di{\NN})$ in a distinct and intrinsic
way by diagram \eqref{eq:intro:LG}.

\medskip
\noindent {\bf 4. Resolutions.}
 The most well-known class of explicit resolutions
of ideals of polynomial rings,
is the Eliahou-Kervaire resolutions of
stable ideals \cite{EK}
(a somewhat more general class than strongly stable ideals),
see also \cite{Pe-St}.

The shift module category enables us to define
a completely new type of projective resolution of strongly stable ideals,
with graded Betti numbers quite distinct
from those in the Eliahou-Kervaire resolution.

The indecomposable projectives in the shift module category are precisely
the strongly stable ideals with a single strongly stable generator.
For instance $x_1x_2^3x_3$ in a polynomial ring $S$ generates a projective
$P = \langle x_1x_2^3x_3 \rangle $.
The ordinary $S$-module minimal free resolution of $P$ is
\[ S(-5)^9 \vpil S(-6)^{12} \vpil S(-7)^4. \]
On the other hand the resolution in the shift module category is simply
$P$, since it is projective. The invariants of these two resolutions
give quite distinct information.
The Betti numbers of the shift module resolution of a strongly stable ideal
reflect more
the combinatorics of the strongly stable generators of the ideal.
As an example class, the shift module resolution of universal lex segment
ideals have Betti numbers like a Koszul resolution, Section \ref{sec:res-ulex}.

\medskip
We also generalize the Eliahou-Kervaire resolution for strongly stable ideals,
to resolutions for a subclass of shift modules,
the {\em rear torsion-free} modules, Section \ref{sec:EK}.


\medskip
The organization of this article is as follows.

\noindent {\bf Part 1.} 
Section \ref{sec:hom-nn}
recalls basic facts on posets, profunctors and the topology on
$\Hom(\NN, \di{\NN})$.

\medskip
\noindent{\bf Part 2.}
Section \ref{sec:st-idealer} gives the correspondence
between strongly stable ideals and open poset ideals in $\Hom(\NN, \di{\NN})$.
We define the duality on strongly stable ideals, and give basic tools
and examples for computing this.
Section \ref{sec:ulex} shows that the duals of universal lex segment ideals
are also universal lex segment.

\medskip
\noindent{\bf Part 3.}
Section \ref{sec:finmod} defines finite shift modules.
Section \ref{sec:finpol} defines shift modules over finite dimensional
polynomial rings. We show that any finitely generated shift module
is derived from a finite shift module by the process of {\em expansion}.
In Section \ref{sec:infinpol} we define shift modules over infinite
dimensional polynomial rings. At the end, in Section \ref{sec:eks} we
give examples of shift modules.

\medskip
\noindent{\bf Part 4.}
Section \ref{sec:finduals} defines duals of finite shift modules,
and Section \ref{sec:polduals} defines duals of shift modules
over polynomial rings. Section
\ref{sec:eksduals} gives examples of duals.

\medskip
\noindent{\bf Part 5.}
Section \ref{sec:eks-res} give examples of how the
ordinary free resolution of a strongly stable ideal and
the shift module resolution differ.
In Section
\ref{sec:res-ulex} we give a shift resolution for strongly stable ideals
which is an analog for the Taylor complex. We give conditions ensuring it
is minimal.
In Section \ref{sec:EK} we establish the Eliahou-Kervaire resolution
for shift modules.

\medskip
\noindent{\bf Appendices.}
Appendix \ref{app:pd-insidens} recalls incidence algebras, and Appendix
\ref{app:ekvi-insidens}
gives the equivalence of categories between shift modules
and modules over the incidence algebras of certain posets.
Appendix \ref{app:closed} states and proves folklore knowledge that
the strongly stable ideals are the most degenerate ideals.

\medskip
\noindent{\it Acknowledgments.} The first steps towards this article
came from a stay of Emily Carrick, McGill University, visiting me at the
University of Bergen. She wrote
programs to compute duals of strongly stable ideals, available
at \cite{EC}, in particular noticing first
the very nice Corollary \ref{cor:st-idealer:xam}. I am grateful
for this initial work, which lead to much more than
first envisioned. 

\part{Isotone maps between natural numbers}
\section{Isotone maps between natural numbers}
\label{sec:hom-nn}
We recall basic notions for partially ordered sets (posets) and
distributive lattices. In particular we consider
isotone maps $f : \NN \pil \di{\NN}$, where $\di{\NN} = \NN \cup \{\infty\}$.
These form themselves a partially ordered set $\Hom(\NN, \di{\NN})$,
which is self-dual. Moreover there is a natural topology on this poset.
We see in Part \ref{part:ss} that it is intimately related to strongly stable
ideals.


\subsection{Posets, isotone maps, and distributive lattices}
\label{subsec:posiso}
Given a poset $P$. A poset ideal $I$ of $P$ is a subset closed
under taking smaller elements, and a poset filter $F$ is a subset closed under
taking larger elements. The distributive lattice $\di{P}$
associated to $P$ is the lattice of all {\it cuts} $(I,F)$
where $I$ is a poset ideal
and $F$ the complement filter of $I$. It is ordered by $(I,F) \leq (J,G)$
if $I \sus J$ (or equivalently $F \supseteq G$). The top element
$(P,\emptyset)$ in $\di{P}$ is denoted $\infty$.

\ignore{
Let $P$ be a partially ordered set, and $P^{\op}$ the opposite
poset. Denote by  $\omega $ the ordered set $\{ d < u \}$.
The distributive lattice associated to $P$ is $\di{P} = \Hom(P^{\op},\omega)$.
The top element of $\di{P}$ is denoted $\ifst$ and sends every
$p^{\op}$ to $u$.

There is a natural map $P \pil \di{P}$ sending an element $p$ to
the element
  \[ \hat{p} : P^{\op} \pil \omega\] where if $p^\prime \geq p$ we send
  $p^{\prime \op} \mapsto d$ and all other elements of $P^\op$ to $u$,
  see Figure \ref{fig:settPoset}.

  \begin{figure}
    \begin{tikzpicture}
 \fill[blue!5] (1,1) ellipse (1 and 2.5);     
\draw (1,1)--(1.8,2.5);
\draw (1,1)--(0.2,2.5);
\draw node[anchor=west] at (1,1) {$p^{\op}$};
\draw node[anchor=west] at (0.7,2.3) {$p^{\prime\op}$};
\draw node at (1.1,1.7) {$d$};
\draw node at (0.8,0.5) {$u$};

\filldraw (1,1) circle (2pt);
\filldraw (0.7,2.3) circle (2pt);
\end{tikzpicture}
\caption{}
\label{fig:settPoset}
\end{figure}
}

  We are particularly interested in this when
  $P$ is totally ordered:
  \[ \NN = \{ 1 < 2 < \cdots \}, \quad [n] = \{ 1 < 2 < \cdots < n \}. \]
Then
\[\di{\NN} = \NN \cup \{ \ifst \}, \quad \diw{[n]} = [n+1] = [n]
  \cup \{ \infty \}. \]

For a poset $P$ denote by $P^\op$ the opposite poset with 
order relation reversed,
so  $p^\op \leq q^\op$ in $P^\op$ if $p \geq q$ in $P$.
For $P,Q$ two posets, a map $f : P \pil Q$ is isotone (order preserving) 
if $p_1 \leq p_2$ implies $f(p_1) \leq f(p_2)$. The set of such maps
is denoted $\Hom(P,Q)$. It is itself a poset by $f \leq g$ if $f(p) \leq g(p)$
for all $p$.
We note that $\di{P}$ naturally identifies with $\Hom(P^\op, \{0 < 1 \})$:
The cut $(I,F)$ corresponds to the morphism
\[ p : P^\op \pil \{ 0 < 1 \}, \quad p^{-1}(0) = F^\op, p^{-1}(1) = I^\op. \]

For $p \in P$ let $\upa p$ be
the filter $\{ q \, | \, q \geq p \}$.
There is a natural isotone map $P \pil \di{P}$ sending
$p \mapsto ((\upa p)^c,\upa p)$. For the totally ordered sets above, this is the
natural inclusion.

    By \cite[Section 2]{Fl-Pro} there is a natural duality
    \[ \Hom(\NN, \di{\NN}) \mto{D} \Hom(\NN, \di{\NN}) \]
    such that $f \leq g$ iff $Df \geq Dg$.
    It makes $\Hom(\NN, \di{\NN})$ into a self-dual poset, i.e. we have 
    an isomorphism:
\begin{equation} \label{eq:settD}
  D : \Hom(\NN, \di{\NN}) \mto{\iso} \Hom(\NN, \di{\NN})^\op.
  \end{equation}
    This duality is easy to explain by an example.

\begin{example} \label{eks:hom-nn:graph}
In Figure \ref{fig:settfDf}
the red  discs form horizontal segments
making the graph of $f$. The values of $f$ are
\[ 2,2,4,5,5,7, \cdots .\]
The blue circles are filled in along vertical segments to make a ``connected
snake'',
and the graph of $Df$ is obtained by considering the vertical axis as the
argument for $Df$ and its graph given by the blue circles. The values of $Df$
are
\[ 1,3,3,4,6,6, \cdots .\]

\begin{figure}
\begin{tikzpicture}
\draw (1,1)--(2,1)--(3,1)--(4,1) --(5,1)--(6,1)--(7,1);
\draw (1,1)--(1,2)--(1,3)--(1,4) --(1,5)--(1,6)--(1,7);
\draw node[anchor=north] at (1,1) {1}
node[anchor=north] at (2,1) {2}
node[anchor=north] at (3,1) {3}
node[anchor=north] at (4,1) {4}
node[anchor=north] at (5,1) {5}
node[anchor=north] at (6,1) {6};

\foreach \x  in {1,...,6}{
        \draw (\x cm, 1cm + 1pt) -- (\x cm,1cm-1pt);
};
        
\draw node[anchor=east] at (1,1) {1}
node[anchor=east] at (1,2) {2}
node[anchor=east] at (1,3) {3}
node[anchor=east] at (1,4) {4}
node[anchor=east] at (1,5) {5}
node[anchor=east] at (1,6) {6};

\foreach \x  in {1,...,6}{
        \draw (1cm + 1pt,\x cm) -- (1cm-1pt,\x cm);
      };
      
\filldraw[color=red] (1,2) circle (2pt);
\filldraw[color=red] (2,2) circle (2pt);
\filldraw[color=red] (3,4) circle (2pt);
\filldraw[color=red] (4,5) circle (2pt);
\filldraw[color=red] (5,5) circle (2pt);
\filldraw[color=red] (6,7) circle (2pt);

\draw[color=blue] (1,1) circle (2pt);
\draw[color=blue] (3,2) circle (2pt);
\draw[color=blue] (3,3) circle (2pt);
\draw[color=blue] (4,4) circle (2pt);
\draw[color=blue] (6,5) circle (2pt);
\draw[color=blue] (6,6) circle (2pt);
\end{tikzpicture}
\caption{}
\label{fig:settfDf}
\end{figure}

\end{example}

\begin{remark} In general isotone maps $P \pil \di{Q}$
  identify as {\it profunctors} $P \promap Q$,
  and are studied in \cite{Fl-Pro}. The duality above
  is a special case of a duality
  \[ D : \text{Prof}(P,Q) \mto{\iso} \text{Prof}(Q,P), \]
  see \cite[Section 2]{Fl-Pro}.
 In the sequel several results from that article are used.
\end{remark}

\subsection{Large and small maps}
This poset of isotone $\NN \pil \di{\NN}$ decomposes into a disjoint union:
\[ \Hom(\NN, \di{\NN}) = \Hom^L(\NN, \di{\NN}) \cup
  \Hom^{u}(\NN, \di{\NN}) \cup \Hom_S(\NN,\di{\NN}), \]
where an isotone map $f : \NN \pil \di{\NN}$ is in
\begin{itemize}
\item $\Hom_S(\NN,  \di{\NN})$ if its values are bounded by a finite
  number in $\NN$. These maps are called {\it small}.
\item $\Hom^L(\NN,\di{\NN})$ if $f(n) = \ifst $ for some $n \in \NN $.
  These maps are called {\it large}.
\item $\Hom^u(\NN,\di{\NN})$ if the image $f(\NN)$ is in
  $\NN = \di{\NN}\backslash \{ \ifst \}$ and is unbounded.
\end{itemize}

The duality \eqref{eq:settD} swaps $\Hom_S(\NN,\di{\NN})$ and
$\Hom^L(\NN, \di{\NN})$ and maps $\Hom^u(\NN,\di{\NN})$ to itself.

\medskip
There is a natural inclusion $[n] \inpil \NN$. It gives
$[n]^\op \inpil \NN^{\op}$. Applying $\Hom(-,\{0 < 1 \})$ we get
$\di{\NN} \pil \diw{[n]}$ where
all $l > n$ are sent to $\ifst \in \diw{[n]}$.

We get a commutative diagram
\begin{equation} \label{eq:settNN}
  \xymatrix{ \Hom(\NN, \di{\NN}) \ar@{>>}[r]  \ar@{>>}[d] & \Hom(\NN, \diw{[n]})
    \ar@{>>}[d] \\
    \Hom([m],\di{\NN}) \ar@{>>}[r] & \Hom([m],\diw{[n]}). }
 \end{equation}
It restricts to diagrams
\begin{equation} \label{eq:hom-nn:diagrams}
  \xymatrix{ \Hom_S(\NN, \di{\NN}) \ar@{>>}[r]  \ar@{>>}[d] & \Hom(\NN, \diw{[n]}) \ar@{>>}[d]\\
    \Hom_S([m],\di{\NN}) \ar@{>>}[r] & \Hom([m],\diw{[n]}) }, \quad
\xymatrix{\Hom^L(\NN, \di{\NN}) \ar@{>>}[r]  \ar@{>>}[d] & \Hom^L(\NN, \diw{[n]}) \ar@{>>}[d] \\
  \Hom([m],\di{\NN}) \ar@{>>}[r] & \Hom([m],\diw{[n]}) }.
\end{equation}
All maps in these diagrams respect duality appropriately,
as is readily seen by Example \ref{eks:hom-nn:graph}.

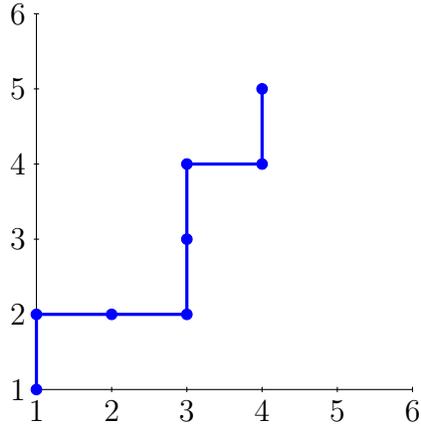
\begin{figure}
\begin{tikzpicture}
\draw (1,1)--(2,1)--(3,1)--(4,1) --(5,1)--(6,1);
\draw (1,1)--(1,2)--(1,3)--(1,4) --(1,5)--(1,6);
\draw node[anchor=north] at (1,1) {1}
node[anchor=north] at (2,1) {2}
node[anchor=north] at (3,1) {3}
node[anchor=north] at (4,1) {4}
node[anchor=north] at (5,1) {5}
node[anchor=north] at (6,1) {6};

\foreach \x  in {1,...,6}{
        \draw (\x cm, 1cm + 1pt) -- (\x cm,1cm-1pt);
};
        
\draw node[anchor=east] at (1,1) {1}
node[anchor=east] at (1,2) {2}
node[anchor=east] at (1,3) {3}
node[anchor=east] at (1,4) {4}
node[anchor=east] at (1,5) {5}
node[anchor=east] at (1,6) {6};

\foreach \x  in {1,...,6}{
        \draw (1cm + 1pt,\x cm) -- (1cm-1pt,\x cm);
      };
      
\filldraw[color=blue] (1,2) circle (2pt);
\filldraw[color=blue] (2,2) circle (2pt);
\filldraw[color=blue] (3,4) circle (2pt);
\filldraw[color=blue] (4,5) circle (2pt);

\filldraw[color=blue] (1,1) circle (2pt);
\filldraw[color=blue] (3,2) circle (2pt);
\filldraw[color=blue] (3,3) circle (2pt);
\filldraw[color=blue] (4,4) circle (2pt);

\draw[very thick, blue] (1,1)--(1,2)--(2,2)--(3,2)--(3,3)--(3,4)--(4,4)--(4,5);
\end{tikzpicture}
\caption{A NE-path from $(1,1)$ to $(4,5)$}
\label{fig:settNE}
\end{figure}

\subsubsection{Large and small maps identify with finite paths}
A finite {\it North-East path} is a path starting form $(1,1)$ going steps of
length one in the north and east directions, see Figure \ref{fig:settNE}.
If it ends at $(m,n)$ we get a function $f : [m] \pil \NN$ with $f(m) = n$,
by letting $f(i)$ be the largest value $j$ such that $(i,j)$ is on the path.
So we have a one-one correspondence
\[ \text{NE-path from } (1,1) \text{ to } (m,n)
  \xleftrightarrow{1-1} \text{isotone maps } f : [m] \pil \NN \text{ with }
  f(m) = n. \]
Let a partial map $f : \NN \dashrightarrow \NN$ be an isotone map
$[m] \pil \NN$, defined for some
initial interval $[m]$ where $m$ is finite. Denote by $\Hom^{P}(\NN,\NN)$
  the set of partial maps. Note that we have a fibration
  \[ \Hom^P(\NN, \NN) \pil \NN^2 \]
  by sending $f$ to $(m,n)$ if $f(m) = n$. The fibers identify as the
  NE paths from $(1,1)$ to $(m,n)$. The cardinality of this
  fiber is $\binom{m+n-2}{n-1}$.

Given a partial map $f$ define two new maps on $\NN$ by
\begin{equation} \label{eq:hom-nn:SL}
 f_S(i) =
  \begin{cases} f(i), & i < m \\
    n+1, & i \geq m
  \end{cases}, \quad
  f^{L}(i) = \begin{cases} f(i), & i \leq m \\
    \ifst, & i > m.
\end{cases}
\end{equation}
  There are then one-to-one correspondences
  \[ \Hom_S(\NN, \di{\NN}) \xleftrightarrow{1-1} \Hom^{P}(\NN, \NN)
    \xleftrightarrow{1-1} \Hom^L(\NN, \di{\NN}). \]

\subsection{The topology on $\Hom(\NN, \di{\NN})$ }
The isotone maps $\Hom(\NN, \di{\NN})$ may be given a topology.
Let $\ul{f}, \ov{f} : \NN \pil \di{\NN}$ be respectively a small
and large isotone map and
\[ U(\ul{f},\ov{f}) = \{ f : \NN \pil \di{\NN}
  \, | \,\, \ul{f} \leq f \leq \ov{f} \}. \]
These sets $U(\ul{f}, \ov{f})$ 
form a basis for a topology on $\Hom(\NN, \di{\NN})$. The map $D$ is
a homeomorphism of topological spaces. In the diagram \eqref{eq:settNN},
we give the other spaces the quotient topology. On $\Hom([m],\diw{[n]})$
this becomes the discrete topology.

\medskip
On a topological space $X$, and $Y$ a subset of $X$, denote by $Y^c$
its complement
$X \backslash  Y $, by $\ov{Y}$ its closure, and by $Y^\circ $
its interior (the union of all open subsets contained in  $Y$).
For any topological space $X$ we have a distinguished subclass of open subsets,
those that are the interiors of their closures. These are called
{\it regular} open sets, and we denote this class as $\regr X$.  
There is an involution, see \cite[Section 4]{Fl-Pro}
\begin{equation} \label{eq:hom-nn-i}
  \regr X \mto{i} \regr X, \quad  U \mapsto (\ov{U})^c = (U^c)^\circ.
  \end{equation}
If $\cI $ is a poset ideal in $\Hom(\NN, \di{\NN})$ both its closure
$\ov{\cI}$ and its interior $(\cI)^\circ $ are poset ideals, and
similarly concerning poset filters, see \cite[Section 4]{Fl-Pro}.
We have the following criterion for
an open poset ideal to be regular.

\begin{definition}
  Let $\cI$ be an open poset ideal and $F: \NN \pil \NN$ an isotone map
  taking finite values. This is a {\it bounding function} for $\cI$ if
  any $f \in \cI$ with $f(p) > F(p)$ is dominated by a $g \in \cI$ (i.e.
  $g \geq f$)  with $g(p) = \infty$. 
  \end{definition}

\begin{proposition} \label{pro:hom-nn:bound}
  An open poset ideal $\cI$ in $\Hom(\NN, \di{\NN})$ is regular
  iff it has a bounding function.
\end{proposition}

\begin{proof}
  By Proposition 6.6b in \cite{Fl-Pro} we have the following criterion:
  Consider sequences  $f_1 \leq f_2 \leq f_3 \leq \cdots$ in $\cI$ and let
  $f$ be $\colim f_i$. Then $\cI$ is regular iff $f$ is in $\cI$
  whenever $f$ is large.

  So suppose $\cI$ has a bounding function. Suppose $f$ large and let
  $m$ such that $f(m-1) < \infty$ and $f(m) = \infty$.
  Since $f_i(j) \leq f(j)$ for
  $j \in [m-1]$ and $f_i(j)$ eventually becomes $f(j)$,
  there is an $N$ such that
  $f_i(j) = f(j)$ for $j \in [m-1]$ and $i \geq N$. We also will have
  $\lim f_i(m) = \infty$ so $f_i(m) > F(m)$ for $i$ large. Then there
  is $g_i \in \cI$ such that $g_i \geq f_i$ and $g_i(m) = \infty$. But
  then $g_i \geq f$ and so $f \in \cI$.

  Conversely assume $\cI$ is regular. Suppose we have defined $F(i)$ for
  $i < m$ such that if $f \in \cI$ and $f(p) > F(p)$ for some $p < m$,
  then there is $g \in \cI$ with $g \geq f$ and $g(p) = \infty$.
  Let \[ T = \{ f \in \cI \, | \, f(m) \text{ finite, and there is no }
    g \in \cI, \, g \geq f \text{ with } g(m) = \infty \}. \]
  If $T = \emptyset$, let $F(m) = F(m-1)$. If $T \neq \emptyset$ consider
  $f \in T$.
  If $f(p) > F(p)$ for some $p < m$, there would be $g \in \cI$ with
  $g \geq f$ and $g(p) = \infty$. This contradicts $f \in T$, so
  $f(p) \leq F(p)$ for every $p < m$.
  
  We show the values of $f(m)$ for $f \in T$ are bounded. Suppose they
  were not.
Then there must be an isotone $\phi : [m-1] \pil \NN$ such that
for any $N$ there is an $f \in T$ with $f(m) \geq N$ and the restriction
$f_{|[m-1]} = \phi$. Let then
\[ f_m(p) = \begin{cases} \phi(p), & p < m \\
    f(m), & p \geq m
  \end{cases}, \quad (f_m \text{ depending on } N). \]
Then clearly $f_m \in T$ since if there is a $g \in \cI$ dominating $f_m$
with $g(m) = \infty$, such a $g$ would contradict $f$ being in $T$.
So we get an increasing sequence of $f_m$'s with 
limit $\phi^L$ (see \eqref{eq:hom-nn:SL}),
which is in $\cI$ since $\cI$ is regular.
But this contradicts $f_m$ being in $T$.

The upshot is that the elements in $T$ have bounded $f(m)$. Let $F(m)$
be the maximal of these.
\end{proof}

\part{Strongly stable ideals and their duality}
\label{part:ss}
\section{Strongly stable ideals and their duals}
\label{sec:st-idealer}

We recall the notion of strongly stable ideals in a polynomial ring
over a field $k$.
For the infinite dimensional polynomial ring $k[x_{\NN}]$ these correspond
precisely to open poset ideals $\cI$ in $\Hom(\NN,\di{\NN})$. Using
the topology on $\Hom(\NN,\di{\NN})$ and that it is a self-dual poset,
we define the dual of
a strongly stable ideal in $k[x_{\NN}]$. We call the ideal dualizable
iff its double dual is the ideal itself. The class of dualizable
strongly stable ideals are those corresponding to regular open poset ideals
$\cI$. We provide results on how to compute the duals of strongly stable
ideals.

\subsection{Correspondence between $\Hom(\NN, \di{\NN})$
  and monomials}
For a set $R$ denote by $k[x_R]$ the polynomial ring in the variables
$\{x_r\}_{r \in R}$ and by $\Mon(x_R)$ the monomials in this ring.
Let $\Mon_{\leq d}(x_R)$ be the monomials of degrees
$\leq d$. There is a bijection \cite[Section 8]{Fl-Pro}
\[ \Hom_S(\NN, \di{\NN}) \mto{\Lambda} \Mon(x_{\NN}), \quad
  f \mapsto \prod_{i \geq 1} x_i^{f(i)-f(i-1)}.\]
Here by convention
$f(0) = 1$. Similarly there is a bijection
\[ \Hom^L(\NN, \di{\NN}) \mto{\Gamma} \Mon(x_{\NN}), \quad
  f \mapsto \prod_{f(i) < \ifst} x_{f(i)}. \]
We get a commutative diagram of bijections, by \cite[Section 8]{Fl-Pro}
\begin{equation} \label{eq:settDLG} \xymatrix{\Hom_S(\NN,\di{\NN}) \ar[dd]^{D} \ar[dr]^{\Lambda} & \\
    & \Mon(x_\NN). \\
    \Hom^L(\NN,\di{\NN}) \ar[ur]^{\Gamma} & 
  }
  \end{equation}
  There are also commutative diagrams of bijections

\begin{equation} \label{eq:st-idealer:mn}
\scalebox{0.8}{ \xymatrix{\Hom(\NN,\diw{[n]}) \ar[dd]^{D} \ar[dr]^{\Lambda} & \\
    & \Mon_{\leq n}(x_\NN) \\
    \Hom([n],\di{\NN}) \ar[ur]^{\Gamma} & 
  } 
\xymatrix{\Hom_S([m],\di{\NN}) \ar[dd]^{D} \ar[dr]^{\Lambda} & \\
    & \Mon(x_{[m]}) \\
    \Hom^L(\NN,\diw{[m]}) \ar[ur]^{\Gamma} & 
  }
\xymatrix{\Hom([m],\di{[n]}) \ar[dd]^{D} \ar[dr]^{\Lambda} & \\
    & \Mon_{\leq n}(x_{[m]}). \\
    \Hom([n],\diw{[m]}) \ar[ur]^{\Gamma} & 
  }  }
\end{equation}

\begin{definition}
If $R$ is a totally ordered set,
a monomial ideal $I$ of the polynomial ring $k[x_R]$ is
{\it strongly stable} (sst) if a monomial $x_j u   \in I$ implies
$x_i u \in I$ for $i < j$. The transitive closure of the relations
i) $x_i u \geq x_ju$ if $i \leq j$ and ii) $x_i u \geq u$ for any $i$
and monomial $u$, gives a partial order $\geq_{st}$ on monomials,
the {\it strongly stable order}. (See also \cite[Lemma 4.2.5]{He-Hi} where
it is called the Borel order.)

If $\{u_i \}_{i \in I}$ is a set
of monomials in $k[x_R]$, we write $\langle u_i \rangle_{i \in I}$ for the
smallest strongly stable monomial ideal containing all the $u_i$. It consists
of all monomials $u$ for which $u \geq_{st} u_i$ for some $i$. We say it is the
{\it strongly stable ideal generated} by the $u_i$'s. We denote by
$\STS(x_{R})$ the strongly stable ideals in $k[x_R]$. (See \cite{FMS}
for more on the perspective of strongly stable generators.)
\end{definition}

The following is immediately verified.
\begin{lemma} \label{lem:st-idealer:order}
  The order relation $\geq_{st}$ corresponds to
  the order relation on $\Hom(\NN, \di{\NN})$ in the following way:
\begin{align} \label{eq:settEkorder}
  f \leq g  \text{ in } \Hom^L(\NN, \di{\NN}) & \ekv \, \Gamma
f \geq_{st} \Gamma g \\ \notag
f \leq g  \text{ in } \Hom_S(\NN, \di{\NN})                                          & \ekv \, \Lambda f \leq_{st} \Lambda g.
\end{align}
\end{lemma}

If $\cI$ is an open poset ideal in $\Hom(\NN, \di{\NN})$, it is fully
determined by its intersection
\[ \cI^L = \cI \cap\Hom^L(\NN, \di{\NN}). \]
Similarly if $\cF$ is an open  poset filter it is fully determined
by
\[\cF_S = \cF \cap \Hom_S(\NN, \di{\NN}).  \]

\begin{proposition} There is a one-to-one correspondence between
  open poset ideals $\cI$ in $\Hom(\NN, \di{\NN})$ and strongly stable
  ideals $I$ in $k[x_\NN]$. For an open poset ideal $\cI $ the associated
  strongly stable ideal is $I = \Gamma(\cI^L)$.
\end{proposition}

\begin{proof}
  By the above observation \eqref{eq:settEkorder}, the image by $\Gamma$ in
  \eqref{eq:settDLG} of $\cI^L$ 
  is a strongly stable ideal $I$ in $\Mon(x_\NN)$.

  When $\cI$ is an open poset ideal, it is fully determined by its elements
  in $\cI^L$, and so distinct $\cI$ give distinct
 strongly stable ideals. Conversely the elements of a strongly stable ideal
 give elements of $\Hom^L(\NN,\di{\NN})$ which generate an open poset ideal
$\cI$ in $\Hom(\NN,\di{\NN})$.
\end{proof}

\begin{remark} \label{rem:st-idealer:cases}
  There are also one-to-one correspondences by $\Gamma$ between
  open poset ideals in $\Hom(\NN,\di{[n]}), \Hom([m],\di{\NN})$ and
  $\Hom([m],\di{[n]]}$ and strongly stable ideals in
  respectively $k[x_{[n]}], k[x_{\NN}]_{\leq m}$ and $k[x_{[n]}]_{\leq m}$.
  \end{remark}
\subsection{Dualizable strongly stable ideals}
If $\cI$ is a {\em regular} open poset ideal of $\Hom(\NN,\di{\NN})$,
let $\cF$
be the poset filter $(\cI^c)^\circ$, the interior of the complement of
$\cI$. Since $i$ of \eqref{eq:hom-nn-i} is an involution, it also
turns $\cF$ into $\cI$.
We say the pair $[\cI,\cF]$ is a {\it Dedekind cut} of $\Hom(\NN, \diw{\NN})$.
Since $D$ is a homeomorphism, the dual $D\cF$ is again a regular open poset
ideal and similarly $D \cI$ a regular open poset filter, and 
$[D \cF, D\cI]$ is the dual Dedekind cut.


\begin{definition} \label{def:st-idealer:dual}
  A strongly stable ideal ideal in $k[x_\NN]$ is
  {\it dualizable strongly stable},
  if it corresponds by $\Gamma$ to a {\it regular} open poset ideal
  $\cI$ in $\Hom(\NN,\di{\NN})$. Its {\it dual strongly stable ideal}
  in $k[y_\NN]$ is the ideal corresponding to the dual regular open poset
  ideal $D\cF$ by the construction above.

  By the commutative diagram \eqref{eq:settDLG} the dual of
  the ideal $\Gamma (\cI^L)$ is $\Lambda(\cF_S)$. 
\end{definition}

\begin{example} \label{eks:st-idealer:dualizable}
  The following two strongly stable ideals are not dualizable.
  \begin{itemize}
  \item $I_1 = \langle x_1, x_2, x_3, \ldots \rangle$,
    the maximal irrelevant ideal in $k[x_\NN]$
\item $I_2 = \langle x_1, x_2^2, x_3^3, x_4^4, \ldots \rangle$
\end{itemize}
 The corresponding open poset ideals $\cI$
are not regular since they do not have bounding functions, Proposition
\ref{pro:hom-nn:bound}.
\begin{itemize}
\item $I_3 = \langle x_1^2, x_1x_2^2, x_1x_2x_3^2, x_1x_2x_3x_4^2, \ldots \rangle.$
  \end{itemize}
  This ideal is dualizable since the corresponding poset ideal is
  regular open. The identity map $\id_\NN$ is a bounding function.
  The ideal $I_3$ is the ideal $I_1^2$ of Example \ref{ex:settGapf}.
  The dual of $I_3$ is $\langle y_1, y_2^2, y_2y_3^2,
  y_2y_3y_4^2, \ldots \rangle$. 
\end{example}

\begin{remark} The above definition is for strongly stable ideals
  in $k[x_\NN]$ but is easily adapted to the more restricted cases.
For instance in $\Hom([m],\di{[n]})$ (which has the discrete topology) let 
$(\cI,\cF)$ a cut. By $\Gamma $ the (open) poset ideal $\cI$ corresponds
to a strongly stable ideal in $k[x_{[n]}]_{\leq m}$.
  The Alexander dual
  poset ideal $\cJ = D \cF$ in $\Hom([n],\di{[m]})$ then gives
  the dual strongly stable ideal $J= \Gamma(\cJ) = \Lambda(\cF)$
  in $k[x_{[m]}]_{\leq n}$. Similarly with other cases of
  Remark \ref{rem:st-idealer:cases}. All maps and correspondences
  in \eqref{eq:hom-nn:diagrams} and \eqref{eq:st-idealer:mn}
  interact well. For instance sst-ideals in $k[x_{[m]}]$ correspond to
  (open) poset ideals in $\Hom(\NN,\diw{[m]})$ (here the topology is the
  discrete topology).
  The dual (open) poset ideal in $\Hom([m], \diw{\NN})$
  then via $\Gamma$ correspond to a strongly stable ideal in
  $k[x_{\NN}]_{\leq m}$. 
\end{remark}

For a Dedekind cut $[\cI, \cF]$ for $\Hom(\NN,\di{\NN})$, its {\it gap} is: 
\[ \cG = \Hom(\NN, \di{\NN}) \backslash (\cI \cup \cF).  \]
The gap is a subset
of $\Hom^u(\NN, \di{\NN})$ by \cite[Cor.4.9]{Fl-Pro}, 
so it consists of unbounded functions.
The following is \cite[Theorem 5.10]{Fl-Pro}.

\begin{theorem} \label{thm:st-id:gap}
  Let $\cI$ be an open poset ideal in $\Hom(\NN, \di{\NN})$.
  Then $\cI$ is also closed iff the strongly
  stable ideal corresponding to $\cI$ is finitely generated.

  So finitely generated strongly stable ideals correspond precisely to
  clopen $\cI$, or alternatively to regular open pairs with empty gap.
  Such strongly stable ideals are therefore dualizable
\end{theorem}

\begin{proof}
By  \cite[Theorem 5.10]{Fl-Pro} $\cI$
is clopen iff it is finitely
generated (and these generators can be chosen large).
But this corresponds to the corresponding strongly stable
ideal being finitely generated.
\end{proof}

Here are examples where the gap consists of a single element.

\begin{example} \label{ex:settGapf}
  Let $\id_{\NN} : \NN \pil \NN$ be the identity function.
  For $1 \leq a < p$ consider the strongly stable ideal $I_a^p$ generated
  by
  \begin{equation} \label{eq:st:iap}
    \{ x_1 x_2 \cdots x_r x_{r+a}^p \, | \, r \geq 0 \}.
    \end{equation}
  Let $J$ be any finitely generated strongly stable ideal, not containing
  the monomial $x_1 x_2 \cdots x_{r-1}x_r$ for any $r \geq 1$.
  Consider the strongly stable ideal $I_a^p + J$ and let $\cI$ be the
  corresponding open poset ideal.
  It is regular by Proposition \ref{pro:hom-nn:bound} and
  Lemma \ref{lem:st-idealer:sums}.
  The monomials of \eqref{eq:st:iap} are
  the $\Gamma h$ of the large isotones $h_r$ defined by
  \[ h_r(i) = \begin{cases} i, & i \leq r \\
      r+a, & r < i \leq r+p \\
      \infty, &  i > r+p
    \end{cases}. \]
  These are then in $\cI$.
  The function $f_r$ below is $\leq h_r$
  \[ f_r(i) = \begin{cases} i, & i \leq r \\
      r+1, & i \geq r+1
    \end{cases} \]
  and so is also in $\cI$. Then $\lim_r f_r = \id_{\NN}$ and is in
  the closure of $\cI$ by \cite[Proposition 5.10]{Fl-Pro}. But $\id_\NN$ is
  not in $\cI$, since if it were, for some $r$ the map
  \[ g_r(i) = \begin{cases} i, & i \leq r \\
      \infty, & i > r
    \end{cases} \]
  would be in $\cI$ due to $\cI$ being open.
  But $\Gamma g_r$ is not in $J$ and not
  in the ideal $I_a^p$ for
  $p > a$. 
 One may argue that $\id_{\NN}$ is the only function in
 the gap for $\cI$.

 Now consider when $p \leq a$.
 Note that when $r = 0$ in \eqref{eq:st:iap}
 the monomial $x_a^p$ is in $I_a^p$.
Since now $p \leq a$ this implies that $x_1 x_2 \cdots x_p$ is in the
  ideal $I_a^p$. All generators of $I_a^p$ with $r \geq p$ is
  a consequence of this, so
  $I_a^p$ in this case is finitely generated.
  \end{example}

  \begin{problem}
    What subsets of $\Hom^u(\NN,\di{\NN})$ can be gaps for Dedekind cuts
    $[\cI,\cF]$?
    \end{problem}

    \subsection{Computing duals of strongly stable ideals}
  
  By Theorem \ref{thm:st-id:gap} any finitely generated ideal is dualizable.
  A package to compute their duals is given in \cite{EC}.
  We describe the duals of principal strongly stable ideals.
  When $A\!: a_1, a_2, \ldots, a_m$ is a finite non-decreasing sequence
  in $\NN$ we get two strongly stable ideals
  \[ \langle y_A \rangle = \langle y_{a_1}y_{a_2} \cdots y_{a_m} \rangle,
    \quad \langle x^A \rangle =
    \langle x_1^{a_1}, x_2^{a_2}, \cdots, x_m^{a_m} \rangle.
    \]

  \begin{proposition} \label{pro:settCI}
    The dual of the strongly stable ideal $\langle y_A \rangle$
    with a single sst-generator, is the strongly stable ideal
    $\langle x^A \rangle$.
  \end{proposition}

  \begin{proof}
    Define the isotone map $f$ by
    \[ f(i) = \begin{cases} a_i, & i = 1, \ldots, m \\
        \ifst, & i > m.
      \end{cases} \]
    Let $\cI$ be the poset ideal generated by $f$, i.e. consisting of all
    isotone maps $g$ such that $g \leq f$. Then $\Gamma(\cI^L)$ is the
    strongly stable ideal generated by
    $\Gamma f = y_{a_1}y_{a_2} \cdots y_{a_m}$.
    The complement filter $\cF$ of $\cI$ is then generated by
    the bounded functions $g_1, \ldots, g_m$ where
    \[ g_p(i) = \begin{cases} 1, & i < p \\
        a_p + 1, & i \geq p.
      \end{cases} \]
    The dual ideal
    is then $J = \Lambda(\cF_S)$, generated by the $\Lambda g_p = x_p^{a_p}$.
    \end{proof}

    \begin{corollary} \label{cor:st-idealer:xam}
      The strongly stable ideals $\langle y_a^m \rangle$
      and $\langle x_m^a \rangle$ are duals of each other.
    \end{corollary}

    \begin{proof} The monomial $y_a^m$ corresponds to the sequence
      $a,a,\cdots, a$ of length $m$. The dual of $\langle y_a^m
      \rangle $ is then the ideal
      $\langle x_1^a, x_2^a, \cdots, x_m^a \rangle$ but this is
      $\langle x_m^a \rangle$.
      \end{proof}

      \begin{remark}
        Strongly stable ideals with one generator are studied in
        \cite{FMS}. In Section 3 they give the minimal
        primary decomposition. In Section 5
        so called Catalan diagrams are introduced associated to
        principal Borel ideals, giving effective computation of Hilbert
        series, and in Section 7 they relate Betti numbers of 
        the principal ideal $\langle x_1x_2\cdots x_n \rangle$
        to pseudo-triangulations.
        In \cite[Example 2]{EK} they note
        that Catalan numbers occur in computing total Betti numbers of
         $\langle x_1, x_2^2, \cdots,  x_n^n \rangle$.
        In \cite{Ar-He} $p$-Borel principal ideals are studied.
        \end{remark}

    The following corollary shows that strongly stable
    duality behaves quite similar to Alexander duality for squarefree
    monomial ideals, see \cite[Cor. 1.5.5]{He-Hi} or \cite[Def.5.20]{Mi-St}.


    \begin{lemma} \label{lem:st-idealer:sums}
      Let $I_1$ and $I_2$ be dualizable strongly stable ideals,
      with duals $J_1$ and $J_2$. Then $I_1 + I_2$ is dualizable strongly
      stable with dual $J_1 \cap J_2$.

      In particular, when $I_1$ is dualizable and $I_2$ is finitely generated,
      then $I_1 + I_2$ is dualizable.
    \end{lemma}

    \begin{proof}
      Let $I_i$ correspond to $\cI_i$, and $J_i$ to $\cJ_i$. Then
      $\cI_1 \cup \cI_2$ is regular open by Proposition \ref{pro:hom-nn:bound}.
      In general the closure $\overline{\cI_1 \cup \cI_2}
      = \overline{\cI_1} \cup \overline{\cI_2}$. Thus the complement
      \[ \overline{\cI_1 \cup \cI_2}^c = \overline{\cI_1}^c \cap
        \overline{\cI_2}^c = \cJ_1 \cap \cJ_2. \]

      In conclusion: $\cI_1 \cup \cI_2$ is regular open with dual regular
      open poset ideal $\cJ_1 \cap \cJ_2$.
    \end{proof}

    \begin{remark} In a topological space $U_1 \cup U_2$ may not
      be regular
      open when $U_1$ and $U_2$ are so. For instance if $[\cI,\cF]$ is a
      regular open pair, then $\cI \cup \cF$ will not be regular open if
      the gap is non-empty.
      \end{remark}


      \begin{corollary} \label{cor:st-idealer:sums}
        Let the $A_1, \ldots, A_r $ each be finite
      weakly increasing sequences of natural numbers.
      The following strongly stable ideals are then duals
      \[ \langle y_{A_1}, \cdots, y_{A_r}
        \rangle, \quad \langle x^{A_1} \rangle \cap \langle x^{A_2} \rangle
        \cap \cdots \cap \langle x^{A_r} \rangle . \]
    \end{corollary}

    In addition to the above, the below seems to be an effective tool
    to compute the dual of a strongly stable ideal. It is
    \cite[Lemma 2.5]{Fl-Pro}.

    \begin{lemma}
      Let $\cI $ be a poset ideal in $\Hom(\NN, \di{\NN})$ and $\cJ$ its
      dual poset ideal. Then $g$ is in $\cJ$ iff for each $f \in \cI$
      there is $p \in \NN$ (depending on $f$)  such that
      $g(f(p)) \leq p$.
    \end{lemma}

    Note: We have no assumptions on these ideals being open or regular.
    
    \begin{proof} By Lemma 2.2 a and b of \cite{Fl-Pro},
      the following holds for
      $p,q \in \NN $ (and is easily checked by Figure \ref{fig:settfDf}):
      \[ q < Dg(p) \text{ iff } g(q) \leq p. \]
      Now $g$ is in $\cJ $ iff $Dg$ is in $\cF$, the complement of $\cI$.
      This holds iff for every $f \in \cI $ there is a $p$ with
      $f(p) < Dg(p)$. By the above this is equivalent to $g(f(p)) \leq p$.
      \end{proof}
  
\section{Universal lex-segment ideals}
\label{sec:ulex}    

If $I \sus k[x_1, \ldots, x_n]$ is a lex segment ideal the extended
ideal $(I) \sus k[x_1, \ldots, x_n, x_{n+1}]$ in a polynomial ring
with one more variable, is usually not lex segment.
However there is a class, universal lex segment ideals, which has
this property. 
They were introduced in \cite{BaNoTh} for finite dimensional polynomial
rings.
We show these ideals correspond precisely to
either an isotone map $f : \NN \pil \NN$ (infinitely generated ideals)
or partial isotones $f: \NN \dashrightarrow \NN$ (finitely generated ideals).
These ideals also have duals, the universal lex segment ideals
corresponding to the dual isotones $Df$.

\medskip
On $\Hom(\NN, \hat{\NN})$ we have the lexicographic order:
$f \ordGL g$ if $f = g$ or $f(r) > g(r)$ where 
\[ r = \min \{n \in \NN \, | \, f(n) \neq g(n)\}. \]
This is a total order which refines the partial order on
$\Hom(\NN, \di{\NN})$. 

\begin{lemma} \label{lem:settDlex}
$f \ordGL g$ iff $Df \ordML Dg$.
\end{lemma}

\begin{proof}
Let $p$ be minimal with $f(p) > g(p)$ and let $m = g(p)$.
Then $Df(m) = p$ and $Dg(m) > p$ as we readily see from
Figure \ref{fig:settfDf}, and $Df(l) = Dg(l)$ for $l < m$.
\end{proof}

Let $f : \NN \pil \NN$ be an unbounded function.
Define a poset ideal and poset filter by (note the order relation is strict)
\begin{align*}
\cI(f) = & \{ g \in \Hom(\NN, \hat{\NN}) \, | \, g \ordM f \}, \\
\cF(f) = & \{ g \in \Hom(\NN, \hat{\NN}) \, | \, g \ordG f \}.
\end{align*}

\begin{proposition}
  $\cI(f)$ and $\cF(f)$ form a Dedekind cut with gap
  $\{f \}$. 
\end{proposition}

\begin{proof}
  The ideal $\cI(f)$ is generated by the large maps
  \begin{equation} \label{eq:ulex:fr}
    f_r (i) = \begin{cases} f(i), & i < r \\
      f(r)-1, & i = r \\
      \infty, & i > r
    \end{cases}, \quad \text{ for } r \geq 1 \text{ and } f(r-1) < f(r), 
    \end{equation}
  and so is an open subset. By Proposition \ref{pro:hom-nn:bound},
$\cI(f)$ is regular. By Lemma \ref{lem:settDlex}
  the dual by $D$ of $\cF(f)$ is $\cI(Df)$. Hence also $\cF(f)$ is regular.
\end{proof}

  \begin{definition} Recall that $\cI^L(f)$ are the large functions
    in $\cI(f)$, and $\cF_S(f)$ the small functions in $\cF(f)$.
    \begin{itemize}
      \item $\tGa(f)$ is the strongly stable ideal $\Gamma(\cI^L(f))$.
\item $\tLa(f)$ is the strongly stable ideal $\Lambda(\cF_S(f))$.
\end{itemize}
\end{definition}

\begin{example}
  In the case when $f = \id_\NN : \NN \pil \NN$ is the identity,
  $\tGa(\id_\NN)$ is the strongly stable ideal $I_1^2$
  (where $a = 1$ and $p = 2$) in Example \ref{ex:settGapf}.

  Note also that
  \[ I_{p-1}^p \supseteq \tGa(\id_\NN) = I_1^2 \supseteq I_1^p. \]
  Thus the $\cI(f)$ are not extremal regular open poset ideals with gap
  $\{f\}$.  
\end{example}

We get the following diagram

\[ \xymatrix{
   & \Hom^u(\NN, \hat{\NN}) \ar[dl]_{\tLa} \ar[dr]^{\tGa} \ar[dd]^D & \\
   \STS(x_\NN) & & \STS(x_\NN) \\
   & \Hom^u(\NN, \hat{\NN}). \ar[ul]^{\tGa} \ar[ur]_{\tLa} & }
\]

\begin{proposition}
The left and right triangles above commute. Furthermore
$\tGa(f)$ and $\tLa(f)$ are dual strongly stable ideals.
\end{proposition}

\begin{proof} By Lemma \ref{lem:settDlex} the dual of
  $\cI(f)$ is $\cF(Df)$. Hence by the commutative
  diagram \eqref{eq:settDLG}, 
  $\Gamma(\cI^L(f)) = \Lambda(\cF_S(Df)$,
  showing that the triangles above commute.

  For the strongly stable ideal $\tGa(f) = \Gamma(\cI^L(f))$, by
  Definition \ref{def:st-idealer:dual} the
  dual ideal is $\Lambda(\cF_S(f)) = \tLa(f)$. 
\end{proof}

Let us describe these ideals more in detail. They may be called (infinitely
generated) universal lex segment ideals, due to the comment after
Proposition \ref{pro:ulex:finulex}.

\begin{proposition} Let the unbounded isotone $f$ take values $f(i) = a_i$, so
$a_1, a_2, a_3, \cdots, $ is the sequence of values.
Recall that we set $a_0 = 1$. 
\begin{itemize}
\item[a.] $\tGa(f)$ is the strongly stable ideal sst-generated by
  \[ x_{a_1} \cdots x_{a_{r-1}}x_{a_r - 1}, \quad \text{for } r \geq 1 \text{ and }
    a_{r-1} < a_r. \] 
\item[b.] Its dual $\tLa(f)$ is the strongly stable ideal sst-generated by
  \[x_1^{a_1 - a_0}x_2^{a_2 - a_1} \cdots x_{r-1}^{a_{r-1} - a_{r-2}} 
x_{r}^{a_r - a_{r-1}+1 }, \quad r \geq 1. \]
\end{itemize}
\end{proposition}

\begin{proof}
  Part a is due to the generators of $\cI(f)$ being \eqref{eq:ulex:fr}.
  Part b follows in a similar way by considering the generators of the
  filter $\cF(f)$. 
\end{proof}


\medskip
The above ideals are not finitely generated since $f$ is unbounded.
Now we consider
universal lex-segment ideal which are finitely generated.
Let $f : [m] \pil \NN$ be an isotone map.
Recall the functions $f_S$ and $f^L$ from \eqref{eq:hom-nn:SL}.
Note that $f_S \ordG f^L$ and this is a covering relation for the
lex order. Define a poset ideal and filter by
\begin{align*}
\cI(f) = & \{ g \in \Hom(\NN, \hat{\NN}) \, | \, g \ordML f^{L} \}, \\
\cF(f) = & \{ g \in \Hom(\NN, \hat{\NN}) \, | \, g \ordGL f_S \}.
\end{align*}

In the same way as above, these form a regular open pair.
In fact a clopen pair.
The associated ideals $\tGa(f)$ and $\tLa(f)$ are dual strongly stable
ideals.

\begin{proposition} \label{pro:ulex:finulex}
  Let the partial isotone $f$ take values
  $a_1, a_2, a_3, \cdots, a_m $ (recall that we set $a_0 = 1$).

 \begin{itemize}
\item[a.] $\tGa(f)$ is the strongly stable ideal generated by
  \[ x_{a_1} \cdots x_{a_{r-1}}x_{a_r - 1},  \quad 1 \leq r \leq m  \text{ and }
  a_{r-1} < a_r \]
together with  $x_{a_1} \cdots x_{a_m}$.

\item[b.] $\tLa(f)$ is the ideal generated by
\[ x_1^{a_1 - a_0}x_2^{a_2 - a_1} \cdots x_{r-1}^{a_{r-1} - a_{r-2}}
x_{r}^{a_r - a_{r-1}+1}, \quad 1 \leq r \leq m. \]
\end{itemize}
These ideals are dual finitely generated universal lex-segment ideals.
\end{proposition}

The description of universal lex segment ideals as above was given
in Proposition 1.2 and Corollary 1.3 in \cite{HiMu-De}, as well
as others characterizations.
The Hilbert functions of these ideals were in \cite{HiMu-Go} characterized
as {\it critical}, i.e. all ideals with this Hilbert functions have the same
Betti numbers.

\part{Shift modules}

We define the category of shift modules. They are the module-theoretic
generalization of strongly stable ideals, much in the same way as
squarefree modules \cite{Ya-Sq} are the generalization of squarefree ideals.
As it turns out, in the natural setting for the duality, our base ring
should be the infinite dimensional polynomial ring. In order to get
there, we go through some steps.

\section{Finite shift modules}

\label{sec:finmod}
Let $\Delta_{m+1}(n)$ be all sequences $\bfd = (d_1, \ldots, d_{m+1})$
of non-negative integers such that the sum $|\dd| = \sum_{i = 1}^{m+1} d_i= n$.
We define shift modules graded by $\Delta_{m+1}(n)$.

\subsection{Finite combinatorial shift modules}

Denote the $i$'th basis vector for ${\mathbb N}^{m+1}$ by $\ee_i$, so
we may write $\bfd$ above as $\sum_{i = 1}^{m+1} d_i\ee_i$. 

\begin{definition}
Let $V$ be a finite dimensional vector space graded by $\Delta_{m+1}(n)$ so
\[ V = \bigoplus_{\dd \in \Delta_{m+1}(n)} V_\dd. \]
This is a
    {\it combinatorial shift module} if
    for each $p = 1, \ldots, m$ and $\dd = (d_1, \ldots, d_{m+1})$  in
    $\Delta_{m+1}(n)$ with $d_{p+1} > 0$ there
    are linear maps \[ s_p : V_{\dd} \pil V_{\dd + \ee_p - \ee_{p+1}} \]
    such that if $1 \leq p < q \leq m$ and $\dd$ has both
    $d_{p+1}, d_{q+1} > 0$,
    the maps $s_p$ and $s_q$ commute:
    \[ s_p \circ s_q = s_q \circ s_p : V_{\dd} \pil V_{\dd + \ee_p + \ee_q -
        \ee_{p+1}- \ee_{q+1}}. \]
If $d_{p+1} = 0$ we define $s_p$ to be zero.

\end{definition}

Homomorphisms $f : V \pil W$ between combinatorial shift modules on
$\Delta_{m+1}(n)$
are then naturally defined by requiring $f$ to commute with the shift maps.
These modules form an abelian category.

For each pair $p < q$ in $[m+1]$ we moreover define shift maps:
\[ s_{p,q} : V_{\dd} \lpil V_{\dd+\ee_p-\ee_q}, \]
as the composition
\[ s_{p,q} = s_{p} \circ s_{p+1} \circ \cdots \circ s_{q-1}. \]

Note then:
\begin{itemize}
\item For $p < r < q$ we have
\[ s_{p,q} = s_{p,r} \circ s_{r,q}. \]
\item When  $d_q$ and $d_\ell$ are both $> 0$, the two maps
  $s_{p,q}$ and $s_{k,\ell}$ commute. In fact they commute save possibly
  in the following cases:
  \begin{itemize}
  \item $d_\ell = 0, d_q > 0$ and $p = \ell$,
    \item $d_q = 0, d_\ell > 0$ and $q = k$.
\end{itemize}
\end{itemize}

\subsection{Relations to the incidence algebra and projectives}
\label{subsec:finmod-proj}
Appendix \ref{app:ekvi-insidens}
shows that for the poset $\Hom([m],\di{[n]})$, the
category of finite dimensional modules over this incidence algebra is
isomorphic to the category of combinatorial shift modules on $\Delta_{m+1}(n)$.
In particular, by Corollary \ref{cor:ekvi:mnpd} this module
category has projective dimension $\min\{m,n\}$.

The projectives in the category of shift modules over $\Delta_{m+1}(n)$
are given as follows: Consider the partial order $\geq_{st}$ on
$\Delta_{m+1}(n)$ given by $\bfd \geq_{st} \bfe$ if
$\sum_{i = 1}^j d_i  \geq \sum_{i = 1}^j e_i$ for $j = 1, \ldots, m+1$.

\begin{lemma}
  $\bfd \geq_{st} \bfe$ iff there is a sequence of shifts taking
  an element of degree $\bfe$ to an element of degree $\bfd$.
\end{lemma}

\begin{proof}
  The if direction is clear since $\bfd + \bfe_p  - \bfe_{p+1}
  \geq_{st} \bfd$. For the only if direction, let $p$ minimal
  such that $d_j > e_j$. Then $\bfd^\prime = \bfd - \bfe_p + \bfe_{p+1}
  \geq_{st} \bfe$, and it follows by induction.
\end{proof}

Then for each $\bfd$ consider the module which is the one-dimensional
vector space $k$ in all degrees
$\bfe \geq_{st} \bfd$ and zero in all other degrees, and with all shift
maps the identity maps. By the equivalence with the incidence algebra,
this is a projective $P(\bfd)$ in the category of combinatorial
shift modules over $\Delta_{m+1}(n)$, and all indecomposable projectives
are of this form.

\subsection{Digression: Algebraic shift modules}
\label{subsec:digression}
Given a combinatorial shift module $V$ on $\Delta_{m+1}(n)$, its associated
{\it algebraic shift module} is defined by the maps, for $p = 1,\ldots, m$:
\begin{equation} \label{eq:finmod:ap}
  a_p:= d_{p+1} \cdot s_p :  V_{\dd} \pil V_{\dd + \ee_p - \ee_{p+1}}
  \end{equation}
where $d_{p+1}$ is the $(p+1)$'th coordinate of
$\dd$. More generally we put $a_{p,q} = d_q \cdot s_{p,q}$.

We then have the commutator relation
\begin{equation} \label{eq:D-acomm} [a_{p,r}, a_{r,q}] = a_{p,q}.
\end{equation}
This makes the module $V$ a module over the Lie algebra $U_{m+1}$ of strictly
upper triangular $(m+1) \times (m+1)$-matrices.

Note however that for any $v \in V_\dd$ the images of the maps included
in relation \eqref{eq:D-acomm} 
\[ a_{p,r} \circ a_{r,q}(v), \quad a_{r,q} \circ a_{p,r}(v), \quad a_{p,q}(v) \]
form at most a one dimensional space, and not two-dimensional (as one would
have expected if one defined algebraic shift maps \eqref{eq:finmod:ap}
as only fulfilling the commutator relation \eqref{eq:D-acomm}).
So the algebraic
shift modules we get from combinatorial shift modules, form a special subclass
of the algebraic shift modules.

In this article we are concerned with combinatorial shifting and
we explain why in Subsection \ref{subsec:finpol-alg}.

\section{Shift modules over the polynomial ring
$k[x_1, x_2, \cdots, x_m]$}
\label{sec:finpol}
We define shift modules over finite dimensional polynomial rings $k[x_{[m]}]$.
When such a module is finitely generated, we show it is induced
by a finite shift module over $\Delta_{m+1}(n)$ for some $n$.
Write $\NN_0^m$ for all $m$-tuples $\dd = (d_1, \ldots , d_m)$ such that
$d_i$ are natural numbers $\geq 0$. For uniformity
of statements we define $d_{m+1} = \infty$ and sometimes consider
an element $\dd$ as $\sum_{i = 1}^{m+1} d_i\bfe_i$. 

\subsection{Combinatorial shift modules}
\begin{definition}
  An $\NN_0^m$-graded vector space $M$ with finite dimensional
  graded parts $M_\dd$ is a {\it combinatorial shift module}
if there for $p = 1, \ldots, m$ are linear maps
\[ s_p : M_\dd \pil M_{\dd + \ee_p - \ee_{p+1}}, \]
whenever $d_{p+1} > 0$, such that $s_p$ and $s_q$ commute
if $d_{p+1}$ and $d_{q+1}$ are both $> 0$. If $d_{p+1} = 0$ we define $s_p$
to be zero.
\end{definition}

For $1 \leq p  < q \leq m+1$ we define
\[s_{p,q} = s_{p} \circ s_{p+1} \circ \cdots \circ s_{q-1}. \]
Again $s_{p,q}$ and $s_{k,\ell}$ commute if
$d_q > 0$ and $d_{\ell} > 0$.

We make a $\NN_0^m$-graded shift module into
a module over the polynomial ring $S = k[x_1, \ldots, x_m]$ as follows.
For an element $u$ of $M_\dd$ define
\[ x_i \cdot u = s_{i,m+1}(u), \quad (\text{recall } d_{m+1} = \infty). \]
Note that for $i < j$:
\[ x_i \cdot u = s_{i,m+1}(u) = s_{i,j} \circ s_{j,m+1} (u)
= s_{i,j} (x_j \cdot u). \]
In particular the polynomial ring $S$ itself becomes a shift module,
by defining
\[ s_p( \cdots x_p^{d_p}x_{p+1}^{d_{p+1}} \cdots ) =
  \cdots x_p^{d_p + 1} x_{p+1}^{d_{p+1} - 1} \cdots. \]
Note that the power $d_{p+1}$ of $x_{p+1}$ does not contribute to
the coefficient of the monomial on the right side, see Subsection
\ref{subsec:finpol-alg}.

The maps $s_p$ are almost $S$-module maps, but not quite.

\begin{lemma} \label{lem:finpol:smod} {}\hskip 1mm
  
  a. Let $u \in M_\dd$ have degree $\dd$. We have 
$s_p(x_i u) = x_i s_p(u)$ except when $i  = p+1$ and $d_{p+1} = 0$. 
In this latter case we have $s_{i-1}(x_i u) = x_{i-1} u$ while
$x_i s_{i-1}(u) = 0$.

b. Generally we have
\begin{equation} \label{eq:finpol:smod}
  s_p(x^{\aa} u) = \begin{cases} 0, & \text{ if } d_{p+1} = 0 \text{ and }
    a_{p+1} = 0, \\
    s_p(x^{\aa}) u, & \text{ if } a_{p+1} > 0,  \\
    x^{\aa} s_p(u), & \text{ if } d_{p+1} > 0. \end{cases}
\end{equation}
Note that when $s_p$ acts non-trivially on both $x^{\aa}$ and $u$
(i.e. both $a_{p+1} > 0$ and $d_{p+1} > 0$), we
can chose which of these to act on:
\[ s_p(x^{\aa} u) = s_p(x^{\aa}) u = x^{\aa} s_p(u). \]
\end{lemma}

\begin{proof}
  a. \begin{eqnarray*} s_p(x_i u) & = & s_p \circ s_{i,m+1}(u) \\
                                  & = & s_{p,p+1} \circ s_{i,m+1}(u)
     \end{eqnarray*}
     If $d_{p+1} > 0$ this is
     \begin{eqnarray*} & = & s_{i,m+1} \circ s_{p,p+1}(u) \\
                       & = & x_i \cdot s_p(u).
     \end{eqnarray*}
     If $d_{p+1} = 0$ and $p+1 \neq i$ then
     \[ s_p(x_i u) = 0 = x_i s_p(u). \]
     If $p+1 = i$ and $d_i = 0$ we have
     \begin{eqnarray*}
       s_{i-1}(x_i u) & = & s_{i-1,i} \circ s_{i,m+1}(u) \\
                      & = & s_{i-1,m+1}(u) \\
                      & = & x_{i-1}u
     \end{eqnarray*}
     and $x_i s_{i-1}(u) = 0$. 

     b. The last property of \eqref{eq:finpol:smod}
     follows due to the $S$-module property in part a. when $d_{p+1} > 0$.
     The middle property of \eqref{eq:finpol:smod} follows by using the
     $S$-module property in part a 
  on $x^{\aa} = x_{p+1}^{a_{p+1}}x^{\aa^\prime}$ as long as $a_{p+1} \geq 2$, and
  then in the last instance using that $s_p(x_{p+1}x^{\aa^\prime}u) =
  x_p x^{\aa^\prime} u$.
\end{proof}

The shift modules over $\NN_0^m$ or equivalently $k[x_1, \ldots, x_m]$
form an abelian category. We denote this category as $\shmod k[x_{[m]}]$.

\medskip


\begin{example}
  Any strongly stable ideal $I$ in $k[x_1, \ldots, x_m]$ is a shift module,
  our primary example of a shift module.
Also any quotient ring $S/I$ of a strongly stable ideal $I$ is a shift
module. More generally if $\{ I_a \}$
is a finite family of strongly stable ideals and
$\{ J_b \}$ is another finite family, 
we get maps 
\[ \oplus_a I_a \pil \oplus_b J_b \]
where the component $I_a \pil J_b$ is either zero or a scalar
multiple of an inclusion map.
The kernel and cokernel of this map are shift modules. We see later
that any shift module is a cokernel of such a map.
\end{example}

A basic result on free shift modules is the following.

\begin{lemma} \label{lem:S-ud}
The free $S$-module $Su_\dd$ with generator $u_\dd$ of degree
$\dd \in \NN_0^m$ can be a shift module iff 
$\dd = d_1 \ee_1$ for some $d_1 \geq 0$.
Then $Su_\dd$ is isomorphic to the strongly stable ideal
$\langle x_1^d \rangle$.
\end{lemma}

\begin{proof} 
  Suppose for some $1 \leq p \leq m-1$ that $d_{p+1} > 0$.
If $Su_\dd$ is a shift module, then note that
$s_p(u_\dd) = 0$ since $Su_\dd$ is zero in degree $\dd + \ee_p - \ee_{p+1}$. 
Thus 
\begin{align*}
  x_p \cdot u_\dd = s_{p,m+1}(u_\dd) & = s_{p,p+1} \circ s_{p+1,m+1}(u_\dd) \\
                                    & = s_{p+1,m+1} \circ s_{p,p+1}(u_\dd)
                                      = s_{p+1,m+1} (0) = 0,
\end{align*}
which is wrong. Here we use that $s_{p,p+1}$ and $s_{p+1,m+1}$
commute since $d_{p+1} > 0$ (and $d_{m+1} = \infty > 0$).

Conversely, if $\dd = d_1 \ee_1$ the map $x^{\aa} u_d \mapsto
x_1^{d_1}x^{\aa}$ defines an isomorphism of shift modules by Lemma
\ref{lem:finpol:smod} b.
\end{proof}

More generally we may prove as above:

\begin{lemma} \label{lem:finpol-Sudd}
  $Su_\dd/(x_1, \ldots, x_{p-1})\cdot u_\dd$ is a shift module
  if $\dd = \sum_{i = 1}^p d_i \ee_i$, but not if
  $d_i > 0$ for some $p+1 \leq i \leq m$.

  This shift module is isomorphic to the cokernel of the inclusion
   \[ \langle x^\dd \cdot x_{p-1} \rangle \inpil \langle x^\dd \rangle. \]
\end{lemma}

\subsection{Projectives}
The following is an immediate consequence of the definition of
shift modules.

\begin{lemma}
  If $s_q \circ s_p(u)$ is non-zero with $q > p$, then
  $s_p \circ s_q(u) = s_q \circ s_p(u)$.
\end{lemma}

\begin{corollary}\label{cor:finpol-qq}
  Given $q_1 \leq q_2 \leq \cdots \leq q_r$ and let $q_1^\prime, \ldots
  q_r^\prime$ be some reordering of these. Then
  \[ s_{q_1} \circ \cdots \circ s_{q_r}(u) =
    s_{q_1^\prime} \circ \cdots \circ s_{q_r^\prime} \]
  if the latter is non-zero.
\end{corollary}

\begin{lemma}\label{lem:finpol-pq}
  If $s_{p_1} \circ \cdots \circ s_{p_r}(u)$ and $s_{q_1} \circ \cdots \circ
  s_{q_t}(u)$ are nonzero of the same degree, with the $p$'s and the
  $q$'s in weakly increasing order, then $r = t$ and each $p_i = q_i$.
\end{lemma}

\begin{proof}
  It is easy to see that we must have $p_1 = q_1$. Then we continue
  by induction.
  \end{proof}

  For a monomial $x^\bfd$ recall that $\langle x^\bfd \rangle$ is
  the strongly stable ideal generated by $x^\bfd$.
  
  \begin{proposition} \label{pro:finpol-pro}
    The ideal $\langle x^\bfd \rangle$ is a projective in the category
    of shift modules over $k[x_{[m]}]$ and all indecomposable
    projectives are of this form.
  \end{proposition}

  \begin{proof}
    Let $M$ be a shift module and $m \in M_{\bfd}$. By Corollary
    \ref{cor:finpol-qq} and Lemma \ref{lem:finpol-pq} there is a unique
    morphism of shift modules $\langle x^\bfd \rangle \pil M$
    sending $x^\bfd \mapsto m$, and
    such that if $x^\bfe = s_{p_1} \circ \cdots \circ s_{p_r}(x^\bfd)$
    then $x^\bfe \mapsto  s_{p_1} \circ \cdots \circ s_{p_r}(m)$.

    Let $M \pil N$ be a surjection, and given $\langle x^\bfd \rangle  \pil N$
    sending $x^\bfd \mapsto n \in N_\bfd$. Let $m \in M_\bfd$ be a lifting
    of $n$. We then get map $\langle x^\bfd \rangle  \pil M$ with
    $x^\bfd \mapsto m$ lifting
    the map $\langle x^\bfd \rangle \pil N$.

    \medskip
    If $P$ is a projective shift module, let $\bfd$ be minimal
    for the order $\geq_{st}$ such that $P_\bfd \neq 0$. Consider
    the short exact sequence
    \[ 0 \pil \langle x^\bfd \rangle_{>_{st}} \pil \langle x^\bfd \rangle
      \pil k \cdot x^\bfd \pil 0. \]
    There is a map $P \pil k\cdot x^\bfd$ which is a map of shift modules.
    It lifts to a map of shift modules $P \pil \langle x^\bfd \rangle$,
    which must be a surjection, and hence the latter is a summand of $P$.
  \end{proof}

\subsection{Digression: Algebraic shifting} \label{subsec:finpol-alg}
For a combinatorial shift module $M$ over $k[x_1, \ldots,x_m]$
its associated {\it algebraic shift maps}
\[ a_p : M_\dd \pil M_{\dd + \ee_p - \ee_{p+1}} \]
are $a_p = d_{p+1} \cdot s_p$ for $p < m$ and
$a_m = s_m$.
More generally define $a_{p,q} = d_q \cdot s_{p,q}$ for $q \leq m$ and
$a_{p,m+1} = s_{p,m+1}$. Again we have the
commutator relation
\[ [a_{p,r} \circ a_{r,q}] = a_{p,q}. \]
Then $M$ becomes a graded module over the Lie algebra
$U_{m+1}$. We also have the natural shift operations $s_p$ and $a_p$
on the polynomial ring $S$. By this the $a_p$ act as derivations:
\[ a_p(x^{\bb} u) = a_p(x^{\bb}) u + x^{\bb} a_p(u). \]


In this article we are only concerned with combinatorial shifting.
The reason is that we want do define dual shift modules
in Part \ref{part:duals}.
A shift map $s_p : M_\dd \pil M_{\dd + \ee_p - \ee_{p+1}}$ dualizes to a shift map
$t_{D_p} : N_\cc \pil N_{\cc + \ee_{D_p } - \ee_{D_p + 1}}$ where
$D_p = 1 + \sum_{i = 1}^p d_i$.
If $d_{p+1} \geq 2$ then $s_p$ and $s_{p+1}$ dualize to
$t_{D_p}$ and $t_{D_{p+1}}$ where the difference $|D_{p+1} - D_p| \geq 2$.
These latter shift maps then commute and so, in a setting where we
have duals,  $s_p$ and $s_{p+1}$ should commute when $d_{p+1} \geq 2$. 
When $d_{p+1} = 1$ we have the possibility
of divergence between combinatorial and algebraic shift modules,
in the presence of duals. 
In this case it is consistent for algebraic shift modules that
$s_p \circ s_{p+1}$ and $s_{p+1} \circ s_p$ are different maps,
as well as their duals $t_{D_p} \circ t_{D_{p+1}}$ and
$t_{D_{p+1}} \circ t_{D_p}$ being different
maps. We stick to combinatorial shift modules where we require also these
two to commute. The category of shift modules is then equivalent
to the category of modules over an incidence algebra,
Appendix \ref{app:ekvi-insidens}.
Since we only consider combinatorial shift modules, we henceforth call
them simply shift modules.

\subsection{Expanded shift modules}
Given a shift module $V$ over $\Delta_{m+1}(n)$ we may extend it to
a shift module $M$ over $k[x_1, \ldots, x_m]$. Let
$\dd = \sum_{i=1}^m d_i \ee_i \in \NN_0^m$.
If the total degree $|\dd| \leq n$, let $d_{m+1} = n - |\dd|$ and
$\hat{\dd} = \dd + d_{m+1} \ee_{m+1}$ which is of degree $n$.
Then let $M_\dd = V_{\hat{\dd}}$.
If $|\dd| \geq n$ then write $\dd = \dd^1 + \dd^2$ where
\begin{equation} \label{eq:finpol:d12}
  \dd^1= (d^1_1, d^1_2, \ldots, d^1_r, 0,\cdots, 0) \quad 
  \dd^2 = (0,\ldots, 0, d^2_r, d^2_{r+1}, \ldots)
  \end{equation}
where the break $r$ is such that
\[ |\dd^1| = n, \quad d_r = d_r^1 + d_r^2, \quad d_r^1 > 0. \]
We let $M_\dd = V_{\dd^1}$. 
Then $M$ becomes a shift module over $\NN_0^m$ as follows.
\begin{itemize}
\item If $|\dd| < n$ let $s^M_p : M_\dd \pil M_{\dd + \ee_p - \ee_{p+1}}$ be
  $s_p : V_{\hat{\dd}} \pil V_{\hat{\dd} + \ee_p - \ee_{p+1}}$. 
\item If $|\dd| \geq n$ let:
  \begin{itemize}
\item If $p < r$ then 
  \[ s_p = s_p^M : M_\dd \pil M_{\dd + \ee_p  - \ee_{p+1}}
    \text{ is } s_p^V : V_{\dd^1} \pil V_{\dd^1 + \ee_p - \ee_{p+1}}. \]
\item If $p \geq r$ then $s_p^M$ is the identity map if $d_{p+1} > 0$,
  otherwise $0$. (Recall when $p = m$ the convention that $d_{m+1} = \infty$.)
\end{itemize}
\end{itemize}

This gives and exact functor from the category of shift
modules over $\Delta_{m+1}(n)$ to the category of shift
modules over $k[x_{[m]}]$. 

\begin{definition}
A shift module $M$ over $k[x_1, \ldots, x_m]$ is {\it expanded}
if it is isomorphic to a module induced from a shift module over
$\Delta_{m+1}(n)$ for some $n$. More specifically it is called
an $n$-expanded shift module.
\end{definition}

We denote the category of $n$-expanded shift modules over $k[x_{[m]}]$ by
$\shmod_{\leq n} k[x_{[m]}]$. So this is a category equivalent to 
the category of shift modules over $\Delta_{m+1}(n)$.
Strongly stable ideals in $k[x_{[m]}]$ generated in degree $\leq n$ are
typical examples of $n$-expanded shift modules.

Let $A,B$ be finite subsets of $\NN_0^m$ and let
\[ \phi : \oplus_{\bfb \in B} \langle x^\bfb \rangle
  \pil \oplus_{\bfa \in A} \langle x^{\bfa} \rangle
\]
be a morphism of shift modules. Note that there is a non-zero morphism
of shift modules $\langle x^\bfb \rangle \pil \langle x^\bfa \rangle$
iff $\bfb \geq_{st} \bfa$ and in this case
$x^\bfb \mapsto \alpha x^\bfb$ for some nonzero $\alpha \in k$.

Let
\[ \hat{b} = \max \{ |\bfb| \, | \, \bfb \in B \}, \quad
  \hat{a} = \max \{ |\bfa| \, | \, \bfa \in A \}. \]

\begin{proposition} Given a morphism $\phi$ which is a minimal
  presentation of its cokernel.
  The regularity of the image $\im \phi$ is $\hat{b}$, and the
  regularity of coker $\phi$ is $\max (\hat{a}, \hat{b}-1)$.
\end{proposition}

\begin{proof}
  Let $n = \hat{b}$. Let $A^\prime \sus A$ be those
  $\bfa \in A$ such that $ |\bfa| \leq n$.
  The image of $\phi$ is contained in $\oplus_{\bfa \in A^\prime}
  \langle x^{\bfa} \rangle$.
  So let $\phi^\prime$ be the map  $\phi$ restricted to this codomain.
  Both $\oplus_{\bfa \in A^\prime} \langle x^{\bfa} \rangle$
  and $\oplus_{\bfb \in B} \langle x^{\bfb} \rangle$ are $n$-expanded,
  and the map $\phi^\prime$ comes from a map
  \[ \overline{\phi^\prime} : \oplus_{\bfb \in B} P(\bfb) \pil
    \oplus_{\bfa \in A^\prime} P(\bfa) \]
  of shift modules over $\Delta_{m+1}(n)$.
  The kernel $\ker \overline{\phi^\prime}$ has finite projective resolution
  with terms being finite sums of projectives shift modules over
  $\Delta_{m+1}(n)$. Each projective $P(\bfc)$ in this resolution must
  have $|\bfc| \leq n$. Since expanding modules is exact, $\im \phi^\prime$
  has a projective resolution $(F_\bullet, d_\bullet)$ of shift modules
  where all terms are projectives $\langle x^\bfc \rangle $ with
  $|\bfc | \leq n$. Such a projective has regularity $|\bfc|$.
  By taking successive mapping cones of this resolution, we get
  that each $\im d_i$ has regularity $\leq n$, and in the end
  $\im \phi^\prime = \im \phi$
  has regularity $\leq n$. But since $\im \phi^\prime$ has a generator
  of degree $\hat{b} = n$, its regularity is $n$.

  Consider the exact sequence
  \[ 0 \pil \im \phi \mto{i} \oplus_{\bfa \in A} \langle x^\bfa \rangle
    \pil \coker \phi \pil 0, \]
  and taking the mapping cone of $i$, the regularity of
  $\coker \phi$ is $\max (\hat{a}, \hat{b} - 1)$.
\end{proof}

\begin{theorem} \label{thm:finpol-exp}
  A shift module $M$ over $k[x_{[m]}]$, finitely generated
  as an $S$-module, is expanded. 
The least $n$ such that $M$ is $n$-expanded
  is either the regularity $\reg \, M$ or $\reg \, M - 1$.
\end{theorem}

\begin{proof} Since $M$ is finitely generated it has a minimal
presentation
  \[  \oplus_{\bfb \in B} \langle x^\bfb \rangle
    \mto{\phi}  \oplus_{\bfa \in A} \langle x^{\bfa} \rangle
\pil M
\]
with $A,B$ finite subsets of $\NN_0^m$. Let $n = \max (\hat{a}, \hat{b})$.
The proof of the previous proposition shows that $\phi$ is induced
from a map of shift modules over $\Delta_{m+1}(n)$:
\[ \oplus_{\bfb \in B} P(\bfb) \mto{\overline{\phi}}
    \oplus_{\bfa \in A} P(\bfa).\]
  Let $\overline{M}$ be the cokernel of $\overline{\phi}$.
  Expanding up to shift modules over $k[x_{[m]}]$ and using that
  expanding is exact, we get that $M$ is expanded from $\overline{M}$. 
  The minimal presentation of $\phi$ is unique up to isomorphism, and
  so we get that the minimal $n$ is $\max(\hat{b}, \hat{a})$.
  If this $n$ is $\hat{a}$, it is the
  regularity of $M$ by the previous proposition.
  If $\hat{a} < \hat{b}$, the regularity
  of $M$ is $\hat{b} - 1 = n-1$.
\end{proof}

\begin{example}
Any $n$-expanded shift module is generated in degrees $\leq n$, but may
be generated in much smaller degrees, for instance if $I$ is a strongly
stable ideal generated in degree $n$, then $S/I$ is generated in degree $0$.
It is $n$-expanded, but
not $(n-1)$-expanded.
\end{example}
  
\section{Shift modules over the infinite dimensional
  polynomial ring $k[x_{\NN}]$}
\label{sec:infinpol}
We define shift modules over the infinite dimensional polynomial
ring $k[x_\NN]$. 
Let $\NN_0^\infty = \oplus_{i \geq 1} \NN_0$ consist of all infinite
sequences $\dd = (d_1, d_2, \ldots, )$ where all $d_i \geq 0$ and
only a finite number of the $d_i$ may be nonzero.
We also put $d_\infty = \infty$.

\begin{definition}
  A $\NN_0^\infty$ graded module $M$ over the infinite dimensional
  polynomial ring is a {\it combinatorial shift
  module} if for every natural number $p \geq  1$ there are
linear maps
\[ s_{p} : M_\dd \pil M_\dd + \ee_p - \ee_{p+1}, \]
with $s_{p}$ the zero map if $d_{p+1} = 0$, such that:
\begin{itemize}
  \item The maps
$s_{p}$ and $s_{q}$ commute if $d_q$ and $d_q$ are both
$> 0$.
\item $s_p(x_{p+1} u) = x_p u$ and if $u$ has degree $\bfd$ with $d_{p+1} > 0$
  this also equals $x_{p+1} s_p(u)$.
  \end{itemize}

For $1 \leq p  < q $ we define
\[s_{p,q} = s_{p} \circ s_{p+1} \circ \cdots \circ s_{q-1}, \]
and $s_{p,\infty}$ to be multiplication by $x_p$. 
Again $s_{p,q}$ and $s_{k,\ell}$ commute if
$d_q > 0$ and $d_{\ell} > 0$.
\end{definition}


\medskip
We could alternatively in a more uniform way
have defined an ${\NN}_0^{\infty}$-graded
shift modules as a graded vector space with maps
for every $1 \leq p < q \leq \infty$
\[ s_{p,q} :  M_\dd \pil M_\dd + \ee_p - \ee_{q}, \]
with $s_{p,q}$ the zero map if $d_q = 0$ such that:
\begin{itemize}
\item $s_{p,r} \circ s_{r,q} = s_{p,q}$, 
\item $s_{p,q}$ and $s_{k,\ell}$ commute when $d_q > 0$ and $d_{\ell} > 0$.
\end{itemize}
In this case $s_{p, \infty}$ would define 
multiplication with $x_p$.

In the same way as before we may also define the algebraic shift maps
$a_{p,q}$ making $M$ an algebraic shift
module, a module over the infinite-dimensional Lie algebra $U_{\infty}$.

\subsection{Extending shift modules}
\medskip
A shift module $M$ over $k[x_1, \ldots, x_m]$ may be extended to
a shift module
\[ \hM = M \te_{k[x_1, \ldots, x_m]} k[x_\NN] \]
  over $k[x_\NN]$. The shift maps are given as:
\[ \hs_{p,\infty}(u \te x^\dd) = \begin{cases} x_p \cdot u \te x^\dd, & p \leq m\\
    u \te x_p \cdot x^\dd, & p > m,  \end{cases} \]
and
\begin{itemize}
  \item For $p < m$
    \[ \hs_p : M_\dd \te x^\cc \mto{s_p \te \en} M_{\dd + \ee_p - \ee_{p+1}} \te x^\cc. \]
\item For $p > m$
\[ \hs_p : M_\dd \te x^\cc\mto{\en \te s_p} M_\dd \te s_p(x^{\cc}). \]
\item For $p = m$, if $x^{\cc}$ contains $x_{m+1}$
  \[ \hs_m : M_\dd \te x^\cc \mto{(\cdot x_m) \te (: x_{m+1})}
    M_{\dd + \ee_m} \te x^{\cc - \ee_{m+1}},\]
  and $\hs_m$ is zero if $x^\cc $
  does not contain $x_{m+1}$.
\end{itemize}

\section{Examples of shift modules}
\label{sec:eks}

In these examples $S$ is a finite dimensional polynomial ring
$k[x_1, \ldots, x_m]$ or infinite dimensional polynomial ring
$k[x_\NN]$. 
\subsection{Strongly stable ideals}

Any strongly stable ideal in a finite or infinite dimensional polynomial
ring is shift module.

Here are examples in the infinite dimensional polynomial ring
$k[x_\NN]$ which we discussed in Example \ref{eks:st-idealer:dualizable}
and will discuss in Section \ref{sec:eksduals}.

\begin{enumerate}
\item The maximal irrelevant ideal in $k[x_\NN]$,
 $I_1 = \langle x_1, x_2, x_3, \ldots \rangle.$
\item $I_2 = \langle x_1, x_2^2, x_3^3, x_4^4, \ldots \rangle.$
\item $I_3 = \langle x_1^2, x_1x_2^2, x_1x_2x_3^2, x_1x_2x_3x_4^2, \ldots \rangle.$
  \end{enumerate}

  \subsubsection{Quotients of strongly stable ideals}
Any quotient module $S/I$ of a strongly stable ideal is a shift module.
More generally if $\{ I_a\}$ is a finite family of strongly stable ideals and
$\{ J_b\}$ is another finite family, 
the quotient by a map
\[ \oplus_a I_a \pil \oplus_b J_b \]
is a shift module, where each component $I_a \pil J_b$ is either
zero or a scalar multiple of an inclusion map.

\subsection{Non-finitely generated shift modules I}
\label{sssec:eksnf1}
For $a \in \NN_0$ let $u_a$ be an element of degree $a \cdot \ee_1$.
\begin{enumerate}
\item 
  The module $M = \oplus_{a \in \NN_0} Su_a$   is a shift module.
\item The module
  \[ M_1 = \oplus_{a \in \NN_0} Su_a/(x_1 u_a)\] is a shift module
  with the same multigraded dimensions as $S$, all are one-dimensional.
  For $M_1$ the shift map $s_1$ is zero.
\end{enumerate}

For $a,b \in \NN_0$ let $u_{a,b}$ be an element of degree
$a \cdot \ee_1 + b \cdot \ee_2$.
\begin{itemize}
\item[(3)] The module
  \[M_2 = \oplus_{(a,b) \in \NN_0^2} Su_{a,b}/(x_1,x_2)u_{a,b} \]
  is again a shift module with the same multigraded dimensions as $S$.
\end{itemize}

\subsection{Non-finitely graded shift modules II}  {} 
\label{sssec:eksnf2} \hskip 1mm

\medskip
\noindent 1. Let the module $N_1$ be the graded vector space which is the direct
sum of all the one-dimensional spaces generated by
$y_\aa = y_{a_1}y_{a_2} \cdots y_{a_r}$ for weakly increasing sequences
$ 1 \leq a_1 \leq a_2 \leq \cdots \leq a_r$, with $r \geq 1$.
As a graded vector space it identifies as $S_{\geq 1}$, the subspace of $S$
spanned by monomials in $S$ of degree $\geq 1$.
We make this a shift
module over $k[y_\NN]$ by
\begin{itemize}
\item[i.] $s_{a_1 - 1} (y_\aa) = 0$
\item[ii.] For $p \geq 2$ if $a_{p-1} < a_p$ then
  \[ s_{a_p-1}(y_\aa) = y_{a_1} \cdots y_{a_p-1} y_{a_{p+1}} \cdots y_{a_r}. \]
\end{itemize}

\medskip
\noindent 2. Let the module $N_2$ be the graded vector space which is the direct
sum of all the one-dimensional spaces generated by
$y_\aa = y_{a_1}y_{a_2} \cdots y_{a_r}$ for weakly increasing sequences
$ 1 \leq a_1 \leq a_2 \leq \cdots \leq a_r$, with $r \geq 2$.
As a graded vector space it identifies as $S_{\geq 2}$, the subspace of $S$
spanned by monomials in $S$ of degree $\geq 2$.
We make this a shift module over $k[y_\NN]$ by:

\begin{itemize}
\item[i.] $s_{a_1 - 1} (y_\aa) = 0$ and $s_{a_2 - 1}(y_\aa) = 0$.
\item[ii.] For $p \geq 3$ if $a_{p-1} < a_p$ then
  \[ s_{a_p-1}(y_\aa) = y_{a_1} \cdots y_{a_p-1} y_{a_{p+1}} \cdots y_{a_r}. \]
\end{itemize}
\part{Duals of shift modules}
\label{part:duals}
For a shift module there is a dual shift module, in much the same way
as for a squarefree module there is a dual squarefree modules. However
shift modules over polynomial  rings are more subtle.
This is due to there not being a simple stable
correspondence between the degrees in the two polynomial rings where
the module and its dual module live over.
    However for finite shift modules this problem is not there and
we first do this case.

\section{Duals of finite shift modules}
\label{sec:finduals}

\subsection{Degree correspondences}

 Recall that $\Delta_{m+1}(n)$ is the set
 of all $(m+1)$-tuples $\bfd = (d_1, \ldots, d_{m+1})$ of integers $\geq 0$
such that their sum is $n$. This set is in bijection with
monomials in $k[x_{[m]}]$ of degree $\leq n$ by sending $\bfd$ to
the monomial $\prod_{i=1}^m x_i^{d_i}$. Recall the rightmost commutative
triangle in \eqref{eq:st-idealer:mn}. Splicing this triangle with its mirror we
get a commutative diagram of bijections.

\begin{equation} \label{eq:dualGLmn} \xymatrix{
  &   \Hom([m],\hat{[n]})  \ar[dl]_{\Gamma} \ar[dd]^D  \ar[dr]^\Lambda &  \\
\Delta_{n+1}(m) \iso \Mon_{\leq m}(y_{[n]})  & & \Mon_{\leq n}(x_{[m]}) \iso
\Delta_{m+1}(n) \\
 &   \Hom([n],\hat{[m]})  \ar[ul]^{\Lambda} \ar[ur]_{\Gamma}  &
 } \end{equation}
The right $\Gamma$-map sends a map $f$ to the monomial
\[  \prod_{i=1}^n x_{f(i)}. \]
For a map $g$ in $\Hom([m],\hat{[n]})$ define $g(0) = 1$
and $g(m+1) = n+1$. Then the right $\Lambda$-map sends a map $g$ to
the monomial
\[ \prod_{i = 1}^m x_i^{g(i) - g(i-1)} \]
and further to the element in $\Delta_{m+1}(n)$
\[ (g(1)-g(0), g(2)-g(1), \cdots, g(m+1) - g(m)). \]

The essential thing for our purpose it that the diagram gives a bijection
between $\Delta_{m+1}(n)$ and $\Delta_{n+1}(m)$. The practical way
to compute this bijection, seems to be using the diagram above. Here is
an example.

\begin{example}
  Consider $(2,1,0,3,1)$ in $\Delta_{4+1}(7)$ corresponding to
  the monomial $x_1^2x_2x_4^3$ in $k[x_{[4]}]_{\leq 7}$.  We want to compute the
  corresponding element in $\Delta_{7+1}(4)$.
The most convenient way is perhaps to use the lower path in \eqref{eq:dualGLmn}.

1. Via $\Gamma $ this monomial $x_1x_1x_2x_4x_4x_4$
corresponds to the function $f$ in
$\Hom([7], \hat{[4]})$ with values $1,1,2,4,4,4,5$. To compute the image
by $\Lambda $ we add $1$ and $5$ at the ends
\[ 1 - 1,1,2,4,4,4,5 - 5 \]
and take the differences $0,0,1,2,0,0,1,0$ giving the element in
$\Delta_{7+1}(4)$, corresponding to the monomial $y_3y_4^2y_7$ in
$k[y_{[7]}]_{\leq 4}$.

\medskip 
We may also use the upper path in $\eqref{eq:dualGLmn}$:

2. From $(2,1,0,3,1)$ we take the partial sums $2,3,3,6$, and add $1$ to each.
Then it corresponds to the function
$g$ in $\Hom([4],\hat{[7]})$ with values $3,4,4,7$. By applying the map
$\Gamma $ this gives $y_3y_4^2y_7$, which in turn may be converted
to the element $(0,0,1,2,0,0,1,0)$ in $\Delta_{7+1}(4)$.
\end{example}

\begin{remark} \label{rem:dualOne}
  The element $1$ in $\Mon_{\leq n}(x_{[m]})$ corresponds
  to $(0,0,\cdots, 0,n)$ in $\Delta_{m+1}(n)$. Via the procedure 1 above
  this corresponds to the function $f$ in $\Hom([n],\di{[m]})$
  with values $m+1, m+1, \cdots, m+1$.
  Taking the differences of
  \[ 1 - (m+1), (m+1), \cdots, (m+1) - (m+1) \]
  we get the sequence $(m,0,0,\cdots, 0)$ giving the monomial
  $y_1^m$ in $k[y_{[n]}]_{\leq m}$.
  In particular we see that as $m$ increases the element $1$ goes to different
  elements.
\end{remark}

\subsection{Duals of finite modules}
For a vector space $V$ denote its dual as $V^* = \Hom(V,k)$.
Let $V = \oplus_{\dd \in \Delta_{m+1}(n)} V_\dd$ be a shift module over
$\Delta_{m+1}(n)$. We want to define a dual shift module
$W = \oplus_{\cc \in \Delta_{n+1}(m)} W_\cc$ over $\Delta_{n+1}(m)$. 
Take $f$ in $\Hom([n],\di{[m]})$ and set
\[ W_{\Lambda f} = (V_{\Gamma f})^*. \]
Next we define the shift maps.
For $p \in [n]$ let $i_p$ be the map with domain $[n]$ such that $i_p(p) = 1$
while all other values of $i_p$ are zero.  Suppose that $f(p) < f(p+1)$
(where we set $f(n+1) = m+1 = \infty$).
Note the shift map
\[ s_p : V_{\Lambda f} \pil V_{\Lambda (f+i_p)}. \]
Let $g = D(f+i_p)$ be the dual. Letting $k = f(p)$ we see by
Figure \ref{fig:dualfpg} that the dual $Df = g + i_k$. 

\begin{figure}
\begin{tikzpicture}
\hskip 1cm \scalebox{0.7}{
\draw (2,1)--(3,1)--(4,1) --(5,1)--(6,1)--(7,1)--(8,1);
\draw (1,2)--(1,3)--(1,4) --(1,5)--(1,6)--(1,7);
\draw[loosely dotted] (1,1)--(2,1);
\draw[loosely dotted] (1,1)--(1,2);
\draw node[anchor=north] at (1,1) {1}
node[anchor=north] at (4,1) {p};

\foreach \x  in {2,...,7}{
        \draw (\x cm, 1cm + 1pt) -- (\x cm,1cm-1pt);
};
\draw (1 cm, 1cm + 3pt)--(1 cm, 1cm - 3pt);     
        
\draw node[anchor=east] at (1,1) {1}
node[anchor=east] at (1,3) {k};

\foreach \x  in {2,...,6}{
        \draw (1cm + 1pt,\x cm) -- (1cm-1pt,\x cm);
      };
\draw (1 cm+3pt, 1cm)--(1 cm-3pt, 1cm);
      
\filldraw[color=red] (2,3) circle (3pt);
\filldraw[color=red] (3,3) circle (3pt);
\filldraw[color=red] (4,3) circle (3pt);
\filldraw[color=red] (5,5) circle (3pt);
\filldraw[color=red] (6,5) circle (3pt);
\filldraw[color=red] (7,7) circle (3pt);

\draw[color=blue] (2,2) circle (3pt);
\draw[color=blue] (5,3) circle (3pt);
\draw[color=blue] (5,4) circle (3pt);
\draw[color=blue] (7,5) circle (3pt);
\draw[color=blue] (7,6) circle (3pt);
}
\end{tikzpicture}
\begin{tikzpicture}
\scalebox{0.7}{
\draw (2,1)--(3,1)--(4,1) --(5,1)--(6,1)--(7,1)--(8,1);
\draw (1,2)--(1,3)--(1,4) --(1,5)--(1,6)--(1,7);
\draw[loosely dotted] (1,1)--(2,1);
\draw[loosely dotted] (1,1)--(1,2);
\draw node[anchor=north] at (1,1) {1}
node[anchor=north] at (4,1) {p};

\foreach \x  in {2,...,7}{
        \draw (\x cm, 1cm + 1pt) -- (\x cm,1cm-1pt);
};
\draw (1 cm, 1cm + 3pt)--(1 cm, 1cm - 3pt);

\draw node[anchor=east] at (1,1) {1}
node[anchor=east] at (1,3) {k};

\foreach \x  in {2,...,6}{
        \draw (1cm + 1pt,\x cm) -- (1cm-1pt,\x cm);
      };
\draw (1 cm+3pt, 1cm)--(1 cm-3pt, 1cm);      
      
\filldraw[color=red] (2,3) circle (3pt);
\filldraw[color=red] (3,3) circle (3pt);
\filldraw[color=red] (4,4) circle (3pt);
\filldraw[color=red] (5,5) circle (3pt);
\filldraw[color=red] (6,5) circle (3pt);
\filldraw[color=red] (7,7) circle (3pt);

\draw[color=blue] (2,2) circle (3pt);
\draw[color=blue] (4,3) circle (3pt);
\draw[color=blue] (5,4) circle (3pt);
\draw[color=blue] (7,5) circle (3pt);
\draw[color=blue] (7,6) circle (3pt);
}
\end{tikzpicture}
\caption{$f$ with ${\color{red} \bullet}$ and $g+i_k$ with
  ${\color{blue} \circ}$. \hskip 1cm
  $f+i_p$ with ${\color{red} \bullet}$ and $g$ with
  ${\color{blue} \circ}$. }
\label{fig:dualfpg}
\end{figure}

Now we have
\[ W_{\Lambda (g+i_k)} = (V_{\Gamma (g+i_k)})^* = (V_{\Lambda D (g+i_k)})^* =
  (V_{\Lambda f})^*. \]
Similarly $W_{\Lambda g} = (V_{\Lambda (f+i_p)})^*$. The dual of the map $s_p$ above
is then by definition our shift map $t_k$ for $W$
\[ t_k : W_{\Lambda g} \pil W_{\Lambda (g + i_k)}. \]

\begin{lemma} \label{lem:dualCD}
  When $V$ is a shift module over $\Delta_{m+1}(n)$,
  the module $W = V^*$ is a shift module over $\Delta_{n+1}(m)$
  with shift maps the $t_k$.
\end{lemma}

\begin{proof} Let $f : [n] \pil \di{[m]}$ be an isotone map. 
  Let $1 \leq p < q \leq n$ with $f(p) < f(p+1)$ and $f(q) < f(q+1)$.
  Since $V$ is a shift module we have a commutative diagram
  \begin{equation} \label{eq:dualVfq}
    \xymatrix{ V_{\Lambda f} \ar[r]^{s_p} \ar[d]^{s_q}
      & V_{\Lambda (f+i_p)} \ar[d]^{s_q} \\
      V_{\Lambda (f + i_q)} \ar[r]^{s_p} & V_{\Lambda (f + i_q + i_p)}.
        }
\end{equation}
      Now let $g = D(f+i_p + i_q)$ and $k = f(p)$ and $\ell  = f(q)$.
      Then as we can infer from Figure \ref{fig:dualfpg}:
      \[ g+ i_\ell = D(f + i_p), \quad g+i_k = D(f+i_q), \quad
        g + i_k + i_\ell = Df. \]
      Dualizing the diagram \eqref{eq:dualVfq} we get a commutative diagram
  \begin{equation*} 
    \xymatrix{ W_{\Lambda (g+i_k + i_\ell)} 
      & W_{\Lambda (g + i_\ell)} \ar[l]^{t_k} \\
      W_{\Lambda (g + i_k)} \ar[u]^{t_\ell}  &
        W_{\Lambda g} \ar[u]^{t_\ell} \ar[l]^{t_k} .
        }
        \end{equation*}
      \end{proof}

      We observe that dualization is an exact functor on the category
      of shift modules over $\Delta_{m+1}(n)$. 

\section{Duals of shift modules over polynomial rings}
\label{sec:polduals}
We now want to define duals of shift modules over
polynomial rings $k[x_{[m]}]$ and $k[x_{\NN}]$.
By the diagram \eqref{eq:dualGLmn}
there is a one-one correspondence between monomials in
$k[x_{[m]}]$ of degree $\leq n$ and monomials in $k[y_{[n]}]$ of degree
$\leq m$. However by Remark \ref{rem:dualOne} such
a correspondence does not stabilize to a correspondence when $n$ and $m$
go to infinity. The element $1$ in the first ring corresponds to the
element $y_1^m$ in the second ring.
Also the diagram \eqref{eq:settDLG} does not extend
to a diamond diagram like \eqref{eq:dualGLmn} since the map $\Lambda$ cannot
be defined on $\Hom^L(\NN, \hat{\NN})$ and the map $\Gamma$ not on
$\Hom_S(\NN, \hat{\NN})$. The notion of dual modules cannot
then rest on a bijection between the monomials in $k[x_{\NN}]$ and $k[y_{\NN}]$.
However we may accomplish our objective using limit considerations.

\medskip
\begin{definition} {} \hskip 1mm
  
  \begin{itemize}
    \item 
  For $f $ in $\Hom^L(\NN, \hat{\NN})$ let $r$ be the integer such that
  $f(r) = \infty $ and $f(r-1) < \infty $. For $m \geq f(r-1)$ define
  $\fm$ to be the small function in $\Hom_S(\NN, \hat{\NN})$ given by
  \[ \fm(i) = \begin{cases} f(i), & i < r \\
      m, & i \geq r
    \end{cases}. \]
  Note that we have 
  \begin{equation} \label{eq:polduals:Lamm}
    \Lambda f^{\underline{m+1}} = x_r \cdot \Lambda \fm.
    \end{equation}
\item   For $g$ in $\Hom_S(\NN, \hat{\NN})$ let $p$ be such that
  $g(p-1) < g(p)$ and $g(i) = g(p)$ for $i \geq p$. For $n \geq p$ let
  $\gn$ be the element in $\Hom^L(\NN, \hat{\NN})$ given by
  \[  \gn (i) = \begin{cases} g(i), & i < n \\
      \infty, & i \geq n
    \end{cases}.
  \]
  Note that we have
  \begin{equation} \label{eq:polduals:Gann}
    \Gamma g_{\mathbf{|}n+1} = x_{g(p)} \cdot \Gamma g_{\mathbf{|}n}.
    \end{equation}
\end{itemize}
\end{definition}

The following is easily verified, since the duality $D$ is essentially
reflection about the axis $x = y$.
\begin{lemma}
  Let $f\in \Hom^L(\NN, \hat{\NN})$ and $g \in \Hom_S(\NN, \hat{\NN})$
  correspond via $D$. Then $\fm \in \Hom_S(\NN, \hat{\NN})$ and
  $g_{\mathbf{|}m} \in \Hom^L(\NN, \hat{\NN})$ correspond via $D$.
\end{lemma}

To define the dual of a module $M$ we need to extend the meaning
of $\Lambda $ to $\Hom^L(\NN, \hat{\NN})$ and of $\Gamma $ to
$\Hom_S(\NN, \hat{\NN})$. 

\begin{definition}
  Let $M$ be a shift module over $k[x_{\NN}]$.
  \begin{itemize}
\item For $g \in \Hom_S(\NN, \hat{\NN})$ we have
  $\gn \in \Hom^L(\NN, \hat{\NN})$. By \eqref{eq:polduals:Gann} we have
multiplication maps
\[ M_{\Gamma g_{\mathbf{|}n}} \mto{\cdot x_{g(p)}} M_{\Gamma g_{\mathbf{|}n+1}}. \]
We define
\[ M_{\Gamma g} = \colim_n M_{\Gamma \gn}.\]
\item  For $f  \in \Hom^L(\NN,\hat{\NN})$ we have
  $\fm \in \Hom_S(\NN, \hat{\NN})$. By \eqref{eq:polduals:Lamm} we have
 multiplication maps
  \[ M_{\ \Lambda \fm} \mto{\cdot x_r} M_{\Lambda  f^{\underline{m+1}}}. \]
  We define
  \[  M_{\Lambda f} = \colim_{m} M_{\Lambda \fm}.\]
\end{itemize}
\end{definition}

\begin{corollary} \label{cor:dualpol}
  For $M$ a shift module over $k[x_\NN]$ and dual elements $f \in Hom^L(\NN, \hat{\NN})$
  and $g \in \Hom_S(\NN, \hat{\NN})$, then
  $M_{\Lambda f} = M_{\Gamma g}$.
\end{corollary}

\begin{proof}
  \[M_{\Lambda f} = \colim_{m} M_{\Lambda \fm} =
      \colim_{m} M_{\Gamma g_{\mathbf{|}m}} = M_{\Gamma g}. \]
\end{proof}

\begin{definition}
  Let $M$ be a shift module over $k[x_{\NN}]$ such that for every
  $g \in \Hom_S(\NN, \di{\NN})$ the graded part $M_{\Gamma g}$ is  
  finite dimensional.
  The {\it dual module} $N = M^\vee$ over $k[y_\NN]$ is defined by letting
  \[ N_{\Lambda g} = (M_{\Gamma g})^* . \]

  Note that if $f = Dg$ is the dual large map, the dual module $N$ also
  by Corollary \ref{cor:dualpol} fulfills
  \[ N_{\Gamma f} = N_{\Lambda g} = (M_{\Gamma g})^* = (M_{\Lambda f})^*.\]
\end{definition}


The shift maps are defined as follows. Let $g \in \Hom_S(\NN, \di{\NN})$.
Suppose $g(p) < g(p+1)$.
For $n$ large enough, we have shift maps
\[ s_{g(p)} : M_{\Gamma (g+ i_p)_{|n}} \pil M_{\Gamma g_{|n}}. \]
Taking the colimit w.r.t. $n$ we get
a shift map and a dual map
\[ s_{g(p)} : M_{\Gamma (g+i_p)} \pil M_{\Gamma g}, \quad
 t_p : N_{\Lambda g} \pil N_{\Lambda (g+i_p)}, \]
where $t_p$ is the shift map for $N$.
That the shift maps $t_p$ and $t_q$ commute when the $(p+1)$'th and $(q+1)$'th
degree coordinates of $\Gamma f$ are positive is checked like in
Lemma \ref{lem:dualCD}.

Since $\colim_{}$ is an exact functor on vector spaces, we observe that
dualization is an exact functor on the category of shift modules over
$k[x_{\NN}]$. 

  \begin{definition}
    A shift module $M$  over $k[x_\NN]$ is {\it dualizable} if it has a dual
    module $N$ such that the dual module of $N$ is $M$ again.
  \end{definition}

  For a shift module $M$ over $k[x_\NN]$ and $u \in \Mon(x_\NN)$ with
  $n \geq \max(u)$ the largest index of a variable in $u$, there
  are maps
  \[ \xymatrix{ M_u \ar[r]^{\cdot x_n^m} \ar[dr]_{x_{n+1}^m} &
      M_{u \cdot x_n^m} \\
      & M_{u \cdot x_{n+1}^m}. \ar[u]_{(s_n)^m}
    } \]
Here $(s_n)^m$ is the shift maps $s_n$ applied $m$ times.
    This gives a diagram
   \[ \xymatrix{ M_u \ar[r] \ar[dr] &  
       \colim_{m} M_{u \cdot x_n^m} \\
      & \colim_{m} M_{u \cdot x_{n+1}^m} \ar[u]_{(s_n)^m}
      } \]   
    and then a map
    \[ M_u \to \lim_{n} \colim_{m}
      M_{u\cdot x_m^n}. \]

    \begin{proposition}
      Let $M$ be a shift module over $k[x_\NN]$ which has a dual module.
      (In particular this holds if there is a uniform bound
      on the dimensions of the $M_u$).
      Then $M$ is dualizable if and only if the
      natural map
      \[ M_u \to {\rm \lim}_{n} {\rm \colim}_{m}
        M_{u\cdot x_n^m} \] is an isomorphism for every $u \in \NN_0^\infty$.
    \end{proposition}

    \begin{proof} 
      Suppose $M$ is dualizable with dual module $N$. Then for large $f$
      \begin{align*}
        M_{\Gamma f} & =  (\colim_{n} N_{\Lambda \fn})^* \\
                     & =  \lim_{n} \, (N_{\Lambda \fn})^*.
      \end{align*}
      But
      \[ (N_{\Lambda \fn})^* = M_{\Gamma \fn} =
        \colim_{m} M_{\Gamma (f^{\underline n})_{|m}}. \]
        If $u = \Gamma f$ then $\Gamma (f^{\underline n})_{|m} = u \cdot x_n^{m-c_0}$
        for suitable $c_0$. Whence
      \begin{align*}
        M_{\Gamma f} & =  \lim_{n} \colim_{m}
                            M_{\Gamma (f^{\underline n})_{|m} } \\
        M_u & =  \lim_{n} \colim_{m}
                  M_{u \cdot x_n^m}.
      \end{align*}
Conversely, suppose the natural map is an isomorphism. Let $N$ be the dual 
module of $M$.
Then
\[ N_{\Lambda f^{\underline n}} = (\colim_m M_{\Gamma (f^{\underline n})_{|m}})^*. \]
  Hence for the dual $M^\prime$ of $N$:
  \[ M^\prime_{\Gamma f} = (N_{\Lambda f})^*
    = \lim_n (N_{\Lambda {f^{\underline n}}})^*
    = \lim_n \,  \colim_m M_{\Gamma (f^{\underline n})_{|m}} = M_{\Gamma f}.
      \]
\end{proof}

\section{Examples of duals}
\label{sec:eksduals}
We give examples of duals of the modules given in Section \ref{sec:eks}.
We first give three examples of ideals where the dual module
is the polynomial ring $S = k[x_{\NN}]$.

\subsection{Duals of ideals}

\noindent {\em $I = S = k[x_\NN]$ {the polynomial ring.}} Then
for $f$ in $\Hom^L(\NN, \hat{\NN})$ with $f(n-1) < f(n) = \ifst$,
the multiplication maps for $m \geq n$
\[ S_{\Lambda \fm} \mto{\cdot x_n} S_{\Lambda f^{\underline{m+1}}} \]
are always isomorphisms of one-dimensional vector spaces.
Hence the dual module of $S$ is $S$ itself.

\medskip
\noindent {\em ${I}$ {is the ideal}
  ${m^r = (x_1, x_2, x_3, \cdots, )^r}$}.
As above the multiplications maps
\[ I_{\Lambda \fm} \mto{\cdot x_n} I_{\Lambda f^{\underline{m+1}}} \]
are always isomorphisms of one-dimensional vector spaces for $m$ sufficiently
large.
Hence the dual module of $I$ is again $S$.

\medskip
\noindent {\em ${ I}$ {is the ideal}
${{ \langle x_1, x_2^2, x_3^3,
\cdots, x_p^p, \cdots \rangle}}.$}
Note that for any monomial $u$ and $m$ large, the maps
\[ I_{u\cdot x_n^m} \mto{\cdot x_n} I_{u\cdot x_n^{m+1}} \]
are isomorphism between one-dimensional spaces for $m \geq n$.
Whence the dual of $I$ is the module $S$.

\medskip
In the last two examples, the ideals were modules which were not dualizable.
When the ideal is dualizable we have the following, which is
analogous to what happens for Alexander duality for squarefree ideals
\cite{Ya-Sq}.

\begin{proposition} \label{pro:eksdIdual}
  Let  $I$ and $J$ be dual strongly stable ideals.
  \begin{itemize}
  \item[a.] The dual of $I$ is the module $S/J$.
  \item[b.] The dual of $S/I$ is the ideal $J$.
  \end{itemize}
\end{proposition}

\begin{proof}
  a.  Let $[\cI,\cF]$ be a Dedekind cut in $\Hom(\NN, \di{\NN})$,
  with $I = \Gamma (\cI^L)$ and $J$ the dual ideal of $I$, which is
  \[ J = \Lambda(\cF_S) = \Gamma( (D\cF)^L).  \]
Note that for $f$ in $\Hom^L(\NN, \di{\NN})$:
\[ (S/J)_{\Gamma f} = \begin{cases} k, & f \not \in (D\cF)^L \\
0, & f \in (D\cF)^L. \end{cases} \]
Let $N$ be the dual shift module of $I$. We show that $N$ lives in
precisely the same degrees as $S/J$ above. Let $f$ an element of
$\Hom^L(\NN, \di{\NN})$ and $g = Df$ the dual small function.
Then
\begin{equation*} N_{\Gamma f} = N_{\Lambda g} = (I_{\Gamma g})^*, \quad
    I_{\Gamma g} = \colim_{n} I_{\Gamma \gn}
    \end{equation*}

\noindent i) Suppose $f$ is in $(D\cF)^L$ so $g$ is in $\cF_S$.
  Note that $\gn$ is then also in $\cF$. Whence $I_{\Gamma \gn} = 0$.
  Then $N_{\Gamma f}$ above is zero.

\medskip
\noindent ii)
If $f$ is not in $(D\cF)^L$ then $g$ is not in $\cF_S$.
Since $g$ is small, it is not
in the gap between $\cI$ and $\cF$, and hence it is in $\cI$.
Since $\cI$ is open, $\gn$ is in $\cI$ for $n$ big, and so
  \[ \colim_n I_{\Gamma \gn} = k. \]
This shows that $N$ and $S/J$ have precisely the same dimensions in each
degree.
  It is readily verified that the shift maps also correspond, so
  $N = S/J$.

  b. The exact sequence
  \[ 0 \pil I \pil S \pil S/I \pil 0 \]
  dualizes to the exact sequence
  \[ 0 \pil (S/I)^\vee \pil S^\vee \pil I^\vee \pil 0. \]
  Here the map between the last  shift modules identify as
  $S \pil S/J$ and so the dual of $S/I$ is $J$.
\end{proof}

\noindent
{\em ${I}$ { the principal strongly stable ideal}
${  \langle  y_{a_1}y_{a_2} \cdots y_{a_r} \rangle}$.} By Proposition
\ref{pro:settCI}
and Proposition \ref{pro:eksdIdual},
the dual module is
$k[x_\NN]/ \langle x_1^{a_1}, x_2^{a_2},\ldots, x_r^{a_r} \rangle$.


  \medskip

\subsection{Duals of modules}

\medskip
\noindent 
{\em The module ${{M = \oplus_{a \in \NN_0} S u_a}}$ } where $u_a$ is a generator
with degree $(a,0,0,\cdots )$. Then $M_{x_1^m}$ is a
vector space of dimension $m+1$ and the maps
\[ M_{x_1^m} \pil M_{x_1^{m+1}} \]
are injections. Thus
\[ \colim_m M_{x_1^{m+1}} \] is not finite-dimensional,
and $M$ has no dual.

\medskip

Recall the modules $N_1$ and $M_1$ in Subsections \ref{sssec:eksnf1} and
\ref{sssec:eksnf2}.

\ignore{
\medskip \noindent
{\em The module ${{M_1 = \oplus_{a \in \NN_0} (S/x_1)\cdot u_a}}.$
Note that $M$ has precisely the same dimensions (dimension $1$)
in all graded pieces as the polynomial ring $S = k[x_\NN]$. The
shift maps $s_p$ are natural for $p \geq 2$ and all the $s_1$ are zero.

Let $N_1$ be the shift module over $k[y_{\NN}]$ which as a vector space
has basis all monomials in $k[y_\NN]$ of non-zero degree, that is all monomials
$y_{a_1}y_{a_2} \cdots y_{a_p}$ where $1 \leq a_1 \leq a_2 \leq \cdots \leq a_p$
with $p \geq 1$. If $a_{r-1} < a_r$, we let the shift map $s_{a_r-1}$
on $N$ be the natural non-zero map for $r \geq 2$. For $r = 1$ we let
$s_{a_1-1}$ be the zero map.
}}

\begin{proposition} The modules $N_1$ and $M_1$ are dual modules.
  Similarly the modules $N_2$ and $M_2$ are dual modules.
  \end{proposition}

  \begin{proof}
    We show that $N_1$ is the dual module of $M_1$.
That $M_1$ 
is the dual module of $N_1$ is an analogous argument.
So let $N^\prime$ be the dual module of $M_1$.
Let $f \in \Hom^L(\NN,\hat{\NN})$. Then
\begin{equation*}
  (N^\prime)_{\Gamma f} = (M_{1, \Lambda f})^*, \end{equation*}
where
\begin{equation} \label{eq:eksMcolim}
  M_{1, \Lambda f} = \colim_{m} M_{1, \Lambda \fm}.
  \end{equation}
If $f$ is the function taking value $\ifst$ at every $i \in\NN$, then
\[ M_{1, \Lambda \fm} = M_{1, x_1^m} \]
and the colimit above is $0$ (the multiplication by $x_1$ on $M_1$ is zero).
Thus $(N^\prime)_{\Gamma f} = (N^\prime)_{{\mathbf 0}} = 0$.

If $f$ is not the above constant function, then for $m$ big then
$\Lambda \fm$ is
$u \cdot x_p^{m-c_0}$ for some fixed monomial $u$ and fixed $p \geq 2$.
Then the colimit in \eqref{eq:eksMcolim} is $k$. Such $f$ correspond to
monomials $u = \Gamma f$ of positive degree and so
$N^\prime_u = k$ for these $u$.

\medskip
As for the shift maps in $N^\prime$
let $f$ have values the finite sequence $a_1,\cdots, a_{r-1},a_{r}, a_{r+1}, \cdots, a_p$,
with $a_{r-1} < a_r$, and let $f^\prime$
have values the finite sequence $a_1,a_2 \cdots , a_{r-1}, a_{r}-1, a_{r+1}, \cdots, a_p$.
The map
\[ t_{a_r - 1} : (N^\prime)_{\Gamma f} \pil (N^\prime)_{\Gamma f^{\prime}} \]
is the dual of the colimit of the shift maps
\[s_r : (M_1)_{\Lambda f^{\prime \underline m}} \pil (M_1)_{\Lambda \fm}. \]
If $r \geq 2$ this is an isomorphism of
one-dimensional spaces. If $r = 1$ this is the zero map.
We find that dual $N^\prime $ of $M_1$ identifies as $N_1$.

In a similar way we show that the dual of $N_1$ is $M_1$.
\end{proof}

\part{Resolutions}

\section{Examples of resolutions}
\label{sec:eks-res}
We give simple examples of minimal projective shift module resolutions,
and in particular see how they differ from the ordinary minimal free
resolutions as $S$-modules. Recall from Proposition \ref{pro:finpol-pro}
that the projective shift modules over the polynomial ring $k[x_{[m]}]$
are the principal strongly stable ideals $\langle x^\bfd \rangle$.

\subsection{Ideals with projective dimension one as shift modules}
Consider the ideal in $k[x_1,x_2,x_3]$ with strongly stable generators
\[  x_1^ax_2^b,\,\,  x_1^{a-1 + r}x_2^{b-r}x_3, \,\, x_1^{a-2 +s}x_2^{b-s}x_3^2, \]
where $r \geq 1$ and $s \geq r+1$. These generators are illustrated
with bullets in Figure \ref{fig:eks-resI}.
The strongly stable ideal generated by these three monomials will have
\begin{itemize}
\item One generator, $x_1^{a+b}$, whose highest index variable is $x_1$,
\item $b$ generators whose highest index variable is $x_2$,
\item $2b+2-r-s$ generators whose highest index variable is $x_3$.
\end{itemize}

The minimal free resolution as $S = k[x_1,x_2,x_3]$-modules,
by Eliahou-Kervaire \cite{EK} has the following form
\begin{equation*} 
  S(-a-b)^{3b+3-r-s} \vpil S(-a-b-1)^{5b+4-2r-2s} \vpil S(-a-b-2)^{2b+2-r-s}.
\end{equation*}

\medskip
On the other hand the minimal projective shift resolution of this ideal
has the following form when
$r \geq 2$ and $s \geq r+2$:
\begin{align*} 
  & \langle x_1^ax_2^b \rangle \oplus \langle x_1^{a-1 + r}x_2^{b-r}x_3 \rangle \oplus
\langle x_1^{a-2 +s}x_2^{b-s}x_3^2 \rangle \\ \notag
\vpil & \langle x_1^{a-1+r}x_2^{b+1-r} \rangle \oplus \langle x_1^{a-2 +s}x_2^{b+1-s}x_3 \rangle
\end{align*}
When  $r = 1$ and $s \geq 3$ the minimal resolution becomes
\begin{equation*} \langle x_1^{a}x_2^{b-1}x_3 \rangle \oplus
\langle x_1^{a-2 +s}x_2^{b-s}x_3^2 \rangle \vpil \langle x_1^{a-2 +s}x_2^{b+1-s}x_3 \rangle.
\end{equation*}

When $r \geq 2$ and $s = r+1 $ the minimal resolution becomes
\begin{equation*} \langle x_1^ax_2^b \rangle \oplus \langle x_1^{a-1 + r}x_2^{b-1-r}x^2_3 \rangle 
\vpil \langle x_1^{a-1 +r}x_2^{b-r}x_3 \rangle.
\end{equation*}

Finally when $r = 1$ and $s=2$ the ideal becomes the projective module
$\langle x_1^{a}x_2^{b-2}x_3^2 \rangle$.

\begin{figure}
\begin{tikzpicture}
\draw (-3,0)--(0,4.8);
\draw (3,0)--(0,4.8) ;
\draw (-1.8,1.92)--(-1.4,1.92);
\draw (-1.4,1.92)--(-0.8,2.88);
\draw (-0.8,2.88)--(-0.4,2.88);
\draw (-0.4,2.88)--(0.4,4.16);

\filldraw[black] (-2.4,0.96) circle (2pt); 
\filldraw[black] (-1.4,1.92) circle (2pt); 
\filldraw[black] (-0.4,2.88) circle (2pt); 

\draw node[anchor=west] at (-3,0){$x_2^{a+b}$};
\draw node[anchor=south] at (0,4.8){$x_1^{a+b}$};
\draw node[anchor=east] at (3,0){$x_3^{a+b}$};
\end{tikzpicture}
\caption{}
\label{fig:eks-resI}
\end{figure}

\medskip
The two resolutions give quite distinct information.
\begin{itemize}
\item By Figure \ref{fig:eks-resI} the Betti numbers in the shift resolution
reflects better the combinatorial nature of the shift generators.
\item For the various pairs $(a,b)$ with $a+b$ fixed, the shift resolutions
above   have the same graded Betti numbers, but those in the EK-resolution vary.
\item It is easy to find strongly stable ideals with the same graded
  Betti numbers in the EK-resolution above (and so in particular have the
  same Hilbert series), but with distinct Betti numbers in the shift
  resolution.
\end{itemize}

In general all strongly stable ideals in three variables
with generators of the same degree, will have shift projective dimension one
or zero.

\subsection{Ideals with projective dimension two
  as shift modules} \label{subsec:pd2}

Consider the ideal with strongly stable generators
\[ x_1^{a+t}x_2^bx_3^c, \,\, x_1^ax_2^{b+r}x_3^c, \,\, x_1^ax_2^bx_3^{c+s}. \]

The minimal projective shift-resolution of this when
$s > r > t$ is 
\begin{align*}
 & \langle x_1^{a+t}x_2^bx_3^c \rangle \oplus \langle x_1^ax_2^{b+r}x_3^c \rangle \oplus
   \langle x_1^ax_2^bx_3^{c+s} \rangle \\
  \vpil &  \langle x_1^{a+t}x_2^{b+r-t}x_3^c \rangle \oplus \langle x_1^ax_2^{b+r}x_3^{c+s-r} \rangle
          \oplus \langle x_1^{a+t}x_2^bx_3^{c+s-t} \rangle \\
  \vpil & \langle x_1^{a+t}x_2^{b+r-t}x_3^{c+s-r} \rangle.
\end{align*}
This is a special case of Proposition \ref{pro:ru:koszul}.

\section{Koszul-type shift resolutions}
\label{sec:res-ulex}
If the sst-generators of a strongly stable ideal are sufficiently generic
we expect the minimal resolution to be given by a Koszul-type resolution:
If there are $n$ generators, the $p$'th Betti number should be
$\binom{n}{p}$. We will see there is always such a resolution for any
set of generators, and give conditions so that it is minimal. In
particular the minimal shift resolution of universal
lex-segment ideals have this form.

\medskip
Recall the strongly stable partial order on $\Mon(x_\NN)$. 
\[ x_{a_1} \cdots x_{a_r} \geq_{st} x_{b_1} \cdots x_{b_s} \]
iff $r \geq s$ and $a_i \leq b_i$ for $i = 1, \cdots, s$.

Recall that $\langle x^\dd \rangle$ and $\langle x^\ee \rangle$
are indecomposable projectives in
$\shmod k[x_\NN]$, where $\dd, \ee$ are in $\NN_0^\infty$.
There is an inclusion of shift modules
$ \langle \dd \rangle \inpil  \langle \ee \rangle$
iff $x^\dd \geq_{\st} x^\ee$.
If the latter does not hold the only map from $\langle \dd \rangle$
to $\langle \ee \rangle$ is the zero map.
By Lemma \ref{lem:st-idealer:order} for $f,g \in \Hom(\NN,\di{\NN})$
there is an inclusion of shift modules
$\langle \Gamma f \rangle \inpil  \langle \Gamma g \rangle$ iff $f \leq g$.

\medskip
\ignore
{For a function $f : \NN \pil \hat{\NN}$ let $J \sus \NN$ be the {\it jump set}
consisting of all $p$ such that $f(p) < f(p+1)$.
For a subset $R$ of the jump set $J$ let $i_R$ be the indicator function
\[ i_R(p) = \begin{cases} 1, & p \in R \\
    0, & p \not \in R.
  \end{cases} \]
Then $f_R = f + i_R$ is also an isotone map from $\NN$ to $\hat{\NN}$.
For $R = \{r_1 < r_2 < \cdots < r_p \}$ let
$\hf_R$ be the small function obtained by making $f_R$ constant
after the maximum $r_p$ in $R$
\[ \hf_R(i ) = \begin{cases} f_R(i), & i \leq r_p \\
 f_R(r_p) , & i \geq r_p. \end{cases} \]}
\ignore{
If $u$ is a monomial, we get a monomial $s_R(u)$ by applying the shift
maps in $S$
\[ s_R (u) = s_{r_1} \circ s_{r_2} \circ \cdots \circ s_{r_p} (u). \]
Note that if $u$ is the monomial $\Lambda f$, then $\Lambda (f + i_R)$
is $s_R(u)$ and $\Lambda(\hf_R)$ is the cutting off of the monomial $s_R(u)$ by removing
all variables $x_q$ where $q > r_p$. }
\ignore{
\medskip
The universal
lex segment ideal $\tLa(f)$ associated to $f$, see Section \ref{sec:ulex},
is the strongly stable ideal generated by all $\Lambda(\hf_{\{ p \}})$ for
$p$ in the jump set $J(f)$ of $f$.
It is finitely generated if $f$ is
a partial map. We give explicitly the minimal free resolution
of the universal lex segment ideal $\tLa(f)$. }

Given isotone maps
\begin{equation} \label{eq:ru:fi}
  f_i : \NN \pil \di{\NN} \quad (\text{or } [m] \pil \di{\NN}),
  \quad i \in J.
\end{equation}
When explicitly mentioned later we might have the following condition
on the $f_i$.

\begin{condition} \label{con:ru:mi}
  For each $i$ there exists $q_i \in \NN$ (or $[m]$) such that
  $f_i(q_i) < f_j(q_i)$ for every $j \in J\backslash \{ i \}$.
  So each $f_i$ has some $q_i$ where $f_i$ is the unique function
  having minimal value at $q_i$.
\end{condition}

Thus $q$ gives an injective function $q : J \pil \NN$ (or $[m]$).
For $R$ a finite subset of $J$ let $f_R$ be the isotone map
which is the meet of the functions $f_i, i \in R$, so:
\[ f_R(p) = \min \{ f_i(p) \, | \, i \in R \}. \]
Let $P(f_R)$ be the projective module $\langle \Gamma(f_R) \rangle$.
There is an inclusion map
  \[ i_{R,S} : P(f_R) \inpil P(f_S)\]
 when $S \sus R$. For $r \in R$ denote by $i_{R,r} =
  i_{R, R \backslash \{r\}}$.

  \medskip
  Let $I$ be the ideal $\langle \Gamma f_i \rangle_{i \in J}$.
  We can now give the terms in the resolution of the quotient ring
  $k[x_\NN]/I$. Let the first term be
  $F_0 = k[x_\NN]$ and for $p \geq 1$ let
  \[ F_p = \oplus_{R \sus J, |R| = p} P(f_R). \]
For $R = \{ r_1 < r_2 < \cdots < r_p \}$, there are  natural maps
\[ P(f_R) \xlongrightarrow{[-i_{R,r_1}, i_{R,r_2}, \cdots, (-1)^p i_{R,r_p} ]}
 \bigoplus_{i = 1}^p P(f_{R\backslash \{ r_i \}}) \sus F_{p-1}. \]
This gives natural maps
\[ F_p \mto{d_p} F_{p-1}, \quad p \geq 1. \]

\begin{proposition} \label{pro:ru:koszul}
  The map $d$ is a differential, and $F_\bullet $ is a
  projective shift resolution of the quotient ring $k[x_\NN]/I$.

  If the $\{f_i\}$ fulfills Condition \ref{con:ru:mi}, then
  $F_\bullet$ is a minimal projective resolution.
\end{proposition}


\begin{remark}
  The above is the analog of the Taylor resolution for square
  free modules \cite[4.3.2, 6.1]{Mi-St}, \cite[Sec.26]{Pe} or
  \cite[Ch.7]{He-Hi}. The condition for minimality
  is also similar to that for the Taylor complex, \cite[Lem.6.4]{Mi-St}.
  \end{remark}

  \begin{proof} 
Let $\sum \alpha_R u_R$ be a syzygy in $F_p$ where the $\alpha_R$ are
nonzero  constants in $k$ and 
 each $\alpha_R u_R$ is in $P(f_R)$. We may assume it is
 homogeneous so all monomials $u_R$ equal a fixed monomial $u$.
 Put a total order on $I$. Order the $p$-subsets of $I$ such
 that $S > R$ if for the maximal $\ell$ such that $s_\ell$ and $r_\ell$ differ,
 we have $s_\ell < r_\ell$. (So $R$ is dragged down by having a heavy rear.)

 Let ${R_0}$ be {\em minimal} among the $R$
 where $\alpha_R$ is nonzero in the sum.
 Write $R_0 = \{r_1 < r_2 < \cdots < r_p \}$. 
Then
\[ P(f_{R_0}) \text{ maps to }
  \bigoplus_{r_i \in R_0} P(f_{R_0 \backslash \{r_i \}}). \]
The image of $\alpha_{R_0}u_{R_0}$
 in $P(f_{R_0 \backslash \{ r_1 \}})$ must cancel against a term in
 the image of $P(f_S)$ for some $S$ occurring in the syzygy $\sum \alpha_Ru_R$.
 So $R_0\backslash \{r_1 \} = S \backslash \{s_t\}$ for some $t$.
 If $t \geq 2$ then $s_i = r_i$ for $i > t$ while $r_t = s_{t-1} < s_t$.
 This contradicts $R_0 < S$. Thus $t = 1$ and since $R_0 < S$ we have
 $r_1 > s_1$. Let $r_0 = s_1$
and $R^\prime = R_0 \cup \{ r_0 \}$. 
If we add plus or minus the image of $d(\alpha_{R} \Gamma(f_{R^\prime}))$
to the syzygy $\sum \alpha_R u_R$, we get a syzygy with larger minimal
$f_R$.
We may thus continue and in the end see that all syzygies
are in the image of  $d$. 

When Condition \ref{con:ru:mi} is fulfilled,
the maps $f_R < f_{R \backslash \{r\}}$ for every $r \in R$. Hence the
  resolution is minimal.
\end{proof}

\begin{corollary} The above gives the minimal shift resolutions of
  the following ideals:
  \begin{itemize}
  \item $\langle x_1, x_2^2, x_3^3, \cdots \rangle $
  \item Universal lex-segment ideals.
  \end{itemize}
\end{corollary}

\section{Generalized Eliahou-Kervaire resolution}
\label{sec:EK}
The resolutions of strongly stable ideals and more generally stable
ideals is the celebrated Eliahou-Kervaire resolution \cite{EK}, a
resolution where the terms and differentials are explicitly described.
See \cite{Pe-St} or \cite[Sec.28]{Pe} for a simple exposition.
Here we generalize this
to shift modules. Another direction where the Eliahou-Kervaire resolution
has recently been generalized is to the resolution of co-letterplace ideals,
\cite{DFN}.

\subsection{Rear torsion-free modules}
For a degree $\dd$ in $\NN_0^\infty$ let $\max(\dd)$ be the largest
index $i$ such that $d_i$ is nonzero. Also let $\min(\dd)$ be the
smallest such index.

\begin{definition}
A shift module $M$ over $k[x_{\NN}]$ is
{\it rear torsion-free} if for every $m \in M_\dd$ and every monomial $x^{\aa}$
with $\max({\dd}) \leq \min({\aa})$, if
$x^{\aa} \cdot m = 0$ then $m = 0$.
\end{definition}

Note that since $M$ is graded by $\NN_0^\infty$ and the $M_\dd$ are
finite-dimensional, the module $M$ has a minimal homogeneous generating set.
Let $\{m_\dd^i\}$ be such a minimal generating set for $M$, with
$m^i_\dd$ of degree $\dd$.

\begin{lemma} Let $M$ be a rear-torsion free module and $m \in M$.
There is a unique way of writing
\[ m = \sum_{i,\dd} \alpha^i_\dd x^{\aa^i_\dd} m^i_\dd  \]
with $\alpha^i_\dd \in k$ and $m^i_\dd \in M_\dd$, 
and for each term $\max(\dd) \leq \min(\aa^i_\dd)$.
\end{lemma}

\begin{proof}
First we do existence. We may in some way
write $m = \sum \alpha^i_\dd x^{\aa^i_\dd}m^i_\dd$.
Consider $x_p \cdot m^i_\dd$ where $p < \max \dd = b$. This is
\begin{align*} s_{p,\infty}(m^i_\dd) & = s_{p,b} \circ s_{b,\infty} (m^i_\dd) \\
 & = s_{b,\infty} \circ s_{p,b}(m^i_\dd) = x_b \cdot s_{p,b}(m^i_\dd).
 \end{align*}
This means that whenever $x_p m^i_\dd$ occurs in a term above, we may replace
it with the term $x_b \cdot s_{p,b}(m^i_\dd)$ where
$b \geq \max{s_{p,b}(m^i_\dd)}$.
Continuing in this way we get the existence of an expression as claimed.

Now consider uniqueness. If we do not have uniqueness, we have a homogeneous
expression of degree $\ee$
\[ 0 = \underset{\aa^i_\dd + \dd= \ee}{\sum} \alpha^i_\dd x^{\aa^i_\dd} m^i_\dd \]
where not all the $\alpha^i_\dd$ are zero.
Let $p = \max(\ee)$. If $x_p$ does not divide $x^{\aa^i_\dd}$ then
we would have $x^{\aa^i_\dd} = 1$, which is not so since the $m^i_\dd$ are
part of a minimal generating set. Hence $x_p$ divides each $x^{\aa^i_\dd}$
By rear torsion-freeness we may divide out by $x_p$ and get
\[  0 = \sum \alpha^i_\dd x^{\aa^i_\dd}/x_p \cdot m^i_\dd. \]
In this way we may continue until some $m^i_\dd$ is a combination of the
other ones, contradicting minimality of the generators.
\end{proof}

\begin{remark}
  Such a unique way of writing an element in the module is
  more or less exactly the same as the unique way of writing
  an element in a quasi-stable submodule in terms of its Pommaret basis,
  \cite[Thm.3.3]{ABRS} or \cite[Prop.4.4, 4.6]{Sei}.

  A difference to Pommaret bases is that those are submodules
  of free modules. In contrast the class of rear torsion-free shift modules
  also includes modules that are not submodules of free modules. For instance
let the multidegree $\bfd = (d_1,d_2,d_3, 0, \ldots, 0)$ have $d_3 > 0$.
Then the module $Su_{\bfd}/(x_1,x_2)u_{\bfd}$ (which is a shift module
by Lemma \ref{lem:finpol-Sudd}) is rear-torsion free if $d_3 > 0$.
However it is not rear torsion-free if $d_3 = 0$.

Another difference is the quasi-stable module are essentially
direct sums of ideals, they are generated by terms $x^{\alpha}e_k$.
In contrast for a shift sub-module of a free module, this may not
be so.
\end{remark}

Let $\{u^i_\dd\}$ be a set of symbols where $u^i_\dd$ has degree $\dd$.
Let $T^i_\dd$ be the subspace of $Su^i_\dd$ with basis
$x^{\aa} u^i_\dd$ where $\max(\dd) \leq \min(\aa)$, and
\[ T = \bigoplus_{i,\dd} T^i_\dd \sus \bigoplus_{i,\dd} S^i_\dd. \]

\begin{corollary} The natural map $T \pil M$ sending $u^i_\dd$ to $m^i_\dd$
is an isomorphism of vector spaces.
\end{corollary}

\begin{proof} This is clear. \end{proof}

For $p \in \NN$ we may then transport the shift map $s_p$ on $M$
to a shift map $s_p$ on $T$. Since $s_p$ and $s_q$ commute on $M_\dd$
when $d_{p+1}$ and $d_{q+1}$ are nonzero, the same holds for $s_p$ and $s_q$
on $T$.
Explicitly we have
\[ s_p(x^\aa \cdot u^i_\dd) = \begin{cases} s_p(x^\aa)  \cdot u^i_\dd, &
p \geq \max(\dd) \\
x^\aa \cdot s_p(u^i_\dd), & p < \max(\dd).
\end{cases}
\]

\subsection{The complex giving the resolution}
Let $F_p$ be the free $S$-module generated by all symbols
$(i_1, \ldots, i_p \, | \, u^i_\dd)$ where
\[ i_1 < i_2 < \cdots < i_p < \max(\dd). \]
This symbol has multidegree $\bfd + \sum_{j = 1}^p e_{i_j}$. 
For a monomial $x^\aa$, we also let $(\ii \, |\, x^{\aa}u^i_\dd)$ be
$x^{\aa} \cdot (\ii \, | \, u^i_\dd)$.
For each $\dd$ choose an arbitrary total order on the $u^i_\dd$'s.
Define a total order on the symbols $(\ii \, | \, u^i_\dd)$ by
$(\jj \, | \, u^j_\dd) > (\ii \, | \, u^i_\ee)$ when
we have the following.
\begin{itemize}
\item If $\dd \neq \ee$ let $p = \max \{ i \, |\, d_i \neq e_i \}$.
  Then we have $d_p > e_p$  (and write also $\dd > \ee$).
  \item If $\dd = \ee$ and $\ii \neq \jj$ let $q = \max \{ r \, | \, i_r \neq j_r \}$.
    Then we have $j_q > i_q$.
    \item If $\dd = \ee$ and $\ii = \jj$ then $u^j_\dd > u^i_\dd$.
    \end{itemize}

In the following $b = \max(\dd)$. 
Define maps $\delta, \mu : F_p \pil F_{p-1}$ by
\begin{align*}
  (i_1,i_2, \cdots, i_p \, | \, u^i_\dd)
  & \overset{\delta}{\mapsto} \sum_q (-1)^q
    x_{i_q} \cdot (i_1, \cdots,\hat{i_q}, \cdots,  i_p \, | \, u^i_\dd) \\
  (i_1,i_2, \cdots, i_p \, | \, u^i_\dd)
  & \overset{\mu}{\mapsto} \sum_q (-1)^q x_b \cdot
    (i_1, \cdots,\hat{i_q}, \cdots,  i_p \, | \, s_{i_q,b}(u^i_\dd)) 
\end{align*}
Note that $s_{i_q,b}(u^i_\dd)$ will typically be rewritten as a linear
combination of products of monomials and other $u^j_\ee$. 

\begin{lemma}
  $d = \delta - \mu $ is a differential, i.e. $d^2 = 0$.
\end{lemma}

This is a simple check using that the maps $s_{i,j}$ commute when they
are nonzero. Note that the free modules $F_p$ are generally not shift
modules, as a free $S$-module only is a shift-module if
its generator has a multidegree $\bfd = d_1e_1$, Lemma \ref{lem:S-ud}.
The following is the generalized Eliahou-Kervaire resolution.

\begin{theorem} Let $M$ be a rear torsion-free shift module.
  The complex $F_\bullet $ is a free resolution $M$.
\end{theorem}

When $M$ is a strongly stable ideal, the resolution $F_\bullet$
is the Eliahou-Kervaire resolution.

\begin{proof}
Given a {\it homogeneous} (for the $\NN_0^\infty$-grading)  syzygy in $F_p$ 
  \begin{equation} \label{eq:EKsum}
    \sum \alpha_{\ii,u} x^{a_{\ii,u}} (\ii \, | \, u).
  \end{equation}
   Let $(\ii^0 \, | \, u^0)$ be maximal of the nonzero terms with respect
  to the order above.

  Note that
  \[ d(\ii^0 \, |\, u^0) = x_{i^0_1} \cdot
    (i^0_2, \cdots, i^0_p \, | \, u^0) + \text{lower terms}. \]
  Since \eqref{eq:EKsum} is a syzygy and so maps by $d$ to zero, in order
  for the above to cancel, we must in \eqref{eq:EKsum} have a term
  $x^{a_{\ii^{0\prime}, u^0}} \cdot (i^{0\prime}_1, i^0_2, \cdots, i^0_p \, | \, u^0)$
  where $i^{0\prime}_1 < i^0_1$. Then we must in \eqref{eq:EKsum}  have
  terms
  \[ \alpha n \cdot x_{i_1^{0\prime}} (\ii^0 \, | \, u^0) -
    \alpha n \cdot x_{i_1^0} (\ii^{0 \prime} \, | \, u^0). \]
  Then subtracting the image of $\alpha n \cdot (i_1^{0 \prime}, i_1^0, i_2^0,
  \cdots, i_p \, | \, u^0)$ from \eqref{eq:EKsum}, we reduce
  \eqref{eq:EKsum} to a syzygy with smaller initial term. We may continue
  until we get zero, and so the kernel of $F_p \mto{d} F_{p-1}$ is
  the image of $F_{p+1} \mto{d} F_p$.
\end{proof}

\begin{remark}
  A similar resolution occurs in W.Seiler \cite[Thm.7.2]{Sei} for the class
  of quasi-stable modules, and generalizing the Eliahou-Kervaire resolution
  for stable ideals, \cite{EK}.
  Again in \cite{Sei} the differential decomposes into two parts,
  which are completely analogous to our $\delta$ and $\mu$.

  However quasi-stable modules are essentially direct sums of ideals.
  So the resolution of \cite[Thm.7.2]{Sei} is essentially
  a resolution of an ideal.

  A difference concerning the terms in the resolution
  is then that our term $S_{i_q,b}(u^i_{\bfd})$ may be
  a linear combination of products of monomials and other basis terms
  $u^i_\bfe$, while in \cite[Thm.7.2]{Sei} the corresponding term
  is only a product of a monomial and a basis element.
  In Section 6 of \cite{Sei} there is
  a more general form for the differential when
  taking resolutions of polynomial submodules. This resolution may however
  not be minimal.
  \end{remark}
\appendix

\section{Incidence algebras}
\label{app:pd-insidens}
Incidence algebras are constructed from partially ordered sets.
They can be viewed as quiver algebras with relations.

Let $P$ be a poset. It gives a quiver with an arrow for each pair
$p^\prime > p$ in $P$ which is a covering relation. Denote this arrow as
$[p^\prime,p]$.
A path
\[  q = p_n > p_{n-1} > \cdots > p_0 = p \]
of covering relations, gives a product
\[ [p_n,p_{n-1}]\cdots [p_2, p_1] \cdot [p_1,p_0] \]
in the quiver algebra.
We form the quiver algebra with relations by setting these products equal for
any two paths from $p$ to $q$. This is the incidence algebra $I(P)$. 

\medskip
A module $M$ over the incidence algebra is a direct sum
$M = \oplus_{p \in P} M_p$ such that for each $q > p$ we have
a map $M_q \vmto{\cdot [q,p]} M_p$ such that all path relations
are respected.
The indecomposable projective modules for the incidence algebra are
the modules, one for each $y \in P$ 
\[ P(y) = \underset{x \geq y}{\prod} k_x \]
where $k_x$ is a copy of $k$ in degree $x$.
The multiplication with $[q,p]$ on $P(y)$
is the identity map from $k_p$ to $k_q$ and is zero on the $k_x$
where $x \neq p$.

\medskip
Let $\hat{P} = \Hom(P^\op,\omega)$ be the associated distributive lattice
to $P$. Then $\hat{P}$ is a Cohen-Macaulay poset by for instance
\cite[Cor.4.5, Ex. 4.6]{BjGaSt}.
By \cite{Po} or \cite{Wo} the incidence algebra $I(\hat{P})$ is then Koszul.
The elements of $\hat{P}$ are poset ideals in $P$ with the ordering
on $\hat{P}$ induced by inclusions of poset ideals.
Let $I \sus J$ be poset ideals with $J\backslash I = \{ x,y \}$ and
$x,y$ incomparable. The ideal of relations for $I(\hat{P})$ are generated
by the quadratic relations as $I$ and $J$ vary
\[ [I \cup \{x,y\},I \cup \{x \}] \cdot [I \cup \{x\}, I] =
[I \cup \{x,y\}, I \cup \{y\}] \cdot [I \cup \{y \}, I]. \]
The Koszul dual $E(\hat{P})$ of $I(\hat{P})$ is then generated by the
relations
\begin{enumerate}
\item For incomparable $x$ and $y$:
\[ [I \cup \{x,y\}, I \cup \{x \}] \cdot [I \cup \{x\}, I] =
- [I \cup \{x,y \}, I \cup \{y\}] \cdot [I \cup \{y \}, I], \]
\item When $y>x$:
  \begin{equation*} \label{eq:pd:ygx}
    [I \cup \{x,y\}, I \cup \{x \}] \cdot [I \cup \{x\}, I] = 0.
    \end{equation*}
\end{enumerate}

\begin{lemma} The largest degree $d$ for which $E(\hat{P})_d$ is nonzero
is the largest cardinality of an antichain in $P$.
\end{lemma}

\begin{proof}
We claim that if  $I \sus J$ are poset ideals and
\[ I = I_0 \sus I_1 \sus \cdots \sus I_n = J \]
a sequence of covering relations (meaning each $I_{p+1}$ has cardinality one
more than $I_p$), then the product
\begin{equation} \label{eq:insidensIJ} [I_n, I_{n-1}] \cdots [I_{1}, I_0]
\end{equation}
is zero iff $J\backslash I$ contains at least two elements $x,y$ which are
comparable by a covering relation, say $y > x$:
By using Relation (1) above repeatedly,
we only change the product \eqref{eq:insidensIJ} by a sign,
and eventually get to a chain where we have successive terms
\[ I_{r+1} = I_r \cup \{ y \}, \quad I_{r} = I_{r-1} \cup \{ x \}. \]
But then by Relation (2) above
this product is zero.
If $J \backslash I$ is an antichain, the product \eqref{eq:insidensIJ} only
can change sign when we take different paths, and is nonzero.
\end{proof}

\begin{corollary} \label{cor:appindGD}
  The global dimension of the incidence algebra $I(\hat{P})$ is the longest
  antichain in $P$.
\end{corollary}

\begin{proof} The semi-simple part $\prod_{p \in\hat{P}} k_p$ of
the incidence algebra $I(\hat{P})$, has minimal
  resolution of length the highest degree in which the Koszul dual algebra
  $E(\hat{P})$ lives. This is due to the resolution of the semi-simple
  part being given by the Koszul dual algebra \cite{BGS},
  and then \cite[I.5.1]{ARS}.
  This degree is the length of the longest antichain in $P$.
\end{proof}

\section{Equivalence with modules over
  incidence algebras} \label{app:ekvi-insidens}
We show that the categories of shift modules are equivalent to
module categories for incidence algebras of the partially ordered
sets that occur in our setting.

\subsection{The finite case} Let $\omega = \{ 0  < 1 \}$.
The distributive lattice $\hat{P}$ then identifies as
$\Hom(P^\op, \omega)$. 
Considering the poset $\Hom([m], \hat{[n]})$ we then have
\[ \Hom([m], \hat{[n]}) = \Hom([m], \Hom([n]^\op,\omega)) = 
  \Hom([m] \times [n]^{\op}, \omega). \]
So this is the distributive lattice associated to $[m] \times [n]^{\op}$.
As a consequence of Corollary \ref{cor:appindGD} we have.

\begin{corollary} \label{cor:ekvi:mnpd}  The global dimension of the
  incidence algebra of
  \begin{itemize}
    \item[a.] $\Hom([m], \hat{[n]})$
  is $\min \{m,n\}$.
  \item[b.] $\Hom(\NN, \hat{[n]})$ and $\Hom([n], \hat{\NN})$ is $n$.
  \end{itemize}
 \end{corollary}

  \begin{proof} a. The longest antichain in $[m]\times [n]^{\op}$ has
    length $\min\{ m,n \}$.
    Similarly the longest antichain in $\NN \times [n]^{\op}$ has length $n$.
  \end{proof}

  Denote the incidence algebra of $\Hom([m], \hat{[n]})$ as
  $I(m,n)$.
  Now given a finite dimensional module
  $M = \oplus_{f \in \Hom([m],\di{[n]})} M_f$
  over this incidence algebra. Recall the map $\Lambda$
  in \eqref{eq:st-idealer:mn} in Section \ref{sec:st-idealer}.
  Let $M_{\Lambda f} = M_f$. 
  We get a vector space $\Lambda M$
  graded by the monomials $\Mon_{\leq n}(x_{[m]})$
    (These monomials are in one-one correspondence
  with $\Delta_{m+1}(n)$.) 
  \[ \Lambda M := \bigoplus_{f \in \Hom([m],\hat{[n]})} M_{\Lambda f}. \]
\begin{proposition} 
    The correspondence $M \pil \Lambda M$ gives an
    isomorphism of categories of finite dimensional modules:
    \[ \text{modules over } I(m,n) \leftrightarrow
      \text{shift modules over }\Delta_{m+1}(n) \iso \shmod_{\leq n}
    k[x_{[m]}].\]
\end{proposition}

  \begin{proof} Let $\alpha$ be such that
    $\alpha(p) < \alpha(p+1)$ and $i_p$ the bump function with takes value
    $1$ at $p$ and zero elsewhere, and $\beta = \alpha + i_p$.
    The multiplication map
    \[  M_\beta \vmto{\cdot [\beta,\alpha]} M_\alpha \]
    in the incidence algebra corresponds to the shift map
    \[ s_p : M_{\Lambda \alpha} \pil M_{\Lambda \beta}. \]
\medskip
If $p < q$ such that $\alpha(q) < \alpha(q+1)$, let
\[ \gamma = \alpha + i_q, \quad \phi = \alpha + i_p + i_q. \]
The relation
\[ [\phi,\beta] \cdot [\beta, \alpha] = [\phi,\gamma] \cdot [\gamma, \alpha] \]
then corresponds to $s_p$ and $s_q$ commuting.
  \end{proof}

  Each element $f$ of the poset $\Hom([m],\di{[n]})$ gives an
  indecomposable projective module $P(f)$ of the incidence algebra
  $I(m,n)$. The module
  in $\shmod_{\leq n} k[x_{[m]}]$ corresponding to $P(f)$ is the principal
  strongly stable ideal $\langle \Lambda f \rangle \sus k[x_{[m]}]$.
  As a consequence of Corollary \ref{cor:ekvi:mnpd} we get:

  \begin{corollary} The global dimension of the module
    category $\shmod_{\leq n} k[x_{[m]}]$ is $\min \{m,n\}$.
  \end{corollary}

  In particular if $I$ is a strongly stable ideal in $k[x_{[m]}]$
  generated by monomials of degree $2$,
  it has projective dimension one in this category (or zero if
  it has only a single strongly stable generator).
  In contrast, in the ordinary category of modules over the
  polynomial ring $S = k[x_{[m]}]$, it may have any projective dimension
  up to $m-1$.
  
  \medskip
  All of the above may be extended to the poset $\Hom([m],\di{\NN})$
  giving an incidence algebra $I(m,\NN)$. 
  Thus the indecomposable projectives in $\shmod k[x_{[m]}]$ are precisely
  the  principal strongly stable ideals for this ring.


  
\medskip We may also use the correspondence $\Gamma$ to get shift modules.
Given again a module $M = \oplus_{f \in \Hom([m],\di{[n]})} M_f$ 
over the incidence algebra $I(m,n)$ we get a shift module over 
$\Delta_{n+1}(m)$:
\[  \Gamma M = \bigoplus_{f \in \Hom([m],\hat{[n]})} (M_{\Gamma f})^*. \]
As above we get:

  \begin{proposition} 
    The correspondence $M \pil \Gamma M$ gives an
    isomorphism of categories of finite dimensional modules
    \[ \text{modules over } I(m,n) \lpil 
      \text{shift modules over }\Delta_{n+1}(m) \iso \shmod_{\leq m} k[x_{[n]}]. \]
  \end{proposition}

  \subsection{Duals}
If $M = \oplus_{p \in P} M_p$ is a module over an incidence algebra
$I(P)$, we get a module $M^\vee$ over the incidence algebra $I(P^{\op})$ of
the opposite poset (where $*$ denotes dual vector space)
\[ M^\vee = \bigoplus_{p^{\op} \in P^{\op}} (M^\vee)_{p^{\op}}
  := \bigoplus_{p \in P} (M_p)^*. \]
Since $\Hom([m],\di{[n]})$ and $\Hom([n],\di{[m]})$ are opposite posets,
by the third diagram of \eqref{eq:st-idealer:mn}, we get a commutative diagram
(modulo identifying 
the double dual $V^{**}$ of a finite dimensional vector space with $V$.) 
\[  \xymatrix{
  &  \mod I(m,n)   \ar[dl]_{\Gamma} \ar[dd]^{()^\vee}  \ar[dr]^\Lambda &  \\
\shmod \Delta_{n+1}(m)   & & \shmod \Delta_{m+1}(n) \\
 &   \mod I(n,m).  \ar[ul]^{\Lambda} \ar[ur]_{\Gamma}  &
} \]

\subsection{The non-finite case}

Consider the poset $\Hom_S(\NN,\hat{\NN})$ of small maps,
and its incidence algebra
$I_S(\NN, \hat{\NN})$. Again we get a functor
\[ \mod I_S(\NN, \hat{\NN}) \mto{\Lambda} \shmod k[x_{\NN}] \]
which is an equivalence of categories.

Moreover for the poset $\Hom^L(\NN,\hat{\NN})$ of large maps
and its incidence algebra
we get an equivalence of categories
\[ \mod I^L(\NN, \hat{\NN}) \mto{\Gamma} \shmod k[x_\NN]. \]
Since $\Hom_S(\NN, \hat{\NN})$ and $\Hom^L(\NN, \hat{\NN})$
are opposite posets, a module $M = \oplus_{f \in \Hom_S(\NN, \di{\NN})} M_{f}$
over $I_S(\NN, \hat{\NN})$ gives
a dual module over $I^L(\NN, \di{\NN})$
\[ M^\vee = \bigoplus_{g \in \Hom^L(\NN, \di{\NN})} (M^\vee)_g
  := \bigoplus_{f \in \Hom_S(\NN,\hat{\NN})} (M_f)^*.\]
where $g = Df$. 

We obtain a commutative diagram
\[  \xymatrix{
   \mod I_S(\NN,\hat{\NN})   \ar[dd]^{()^*}  \ar[dr]^\Lambda &  \\
 & \shmod k[x_\NN] \\
  \mod I^L(\NN,\hat{\NN}).  \ar[ur]_{\Gamma}  &
} \]

\section{The most degenerate ideals}
\label{app:closed}
In the introduction we stated that the strongly stable ideals
are the most degenerate ideals in a polynomial ring (characteristic $k$ is $0$).
We state this in precise form and give the argument as it seems not easy to
come by in the literature. 

The group $GL(n+1)$ of invertible linear operators on the linear space
generated by the variables, acts by coordinate change on homogeneous ideals in
$I \sus k[x_0, \ldots, x_n]$, where in this appendix $k$ may have any
characteristic. Let $B = B(n+1)$ be the Borel subgroup of upper
triangular matrices of $G = GL(n+1)$, those invertible linear maps sending
$x_j \mapsto \sum_{i = 1}^j \alpha_{ij} x_i$, where the $\alpha_{ij} \in k$.
An ideal $I$ is {\it Borel-fixed} if $g.I = I$ for every $g \in B$.
When char. $k = 0$ this is the same as $I$ being strongly stable,
\cite[Prop. 4.2.4]{He-Hi}. 

Let $\Hilb^{\Pn}$ be the Hilbert scheme of subschemes of the projective
space $\Pn$. We get the action
\begin{equation} \label{app:degen}
  GL(n+1) \times \Hilb^{\Pn} \pil \Hilb^{\Pn}
\end{equation}
where if $x$ corresponds to the ideal $I$ then $g.x$ corresponds to $g.I$.

\begin{theorem} \label{thm:degen} {} \hskip 8mm 

  {\rm a)} The closed orbits of the action \eqref{app:degen} are precisely
  the orbits of Borel-fixed ideals.

  {\rm b)} (char.$\,k = 0$) Any such orbit has exactly one Borel-fixed ideal.

\end{theorem}

\noindent{\it Note.} Part b) is likely true in arbitrary characteristic
but needs a more elaborate proof.

\begin{proof}
  Part a). Given $x \in \Hilb^{\Pn}$
  we get a morphism
  \[ G \lpil \Hilb^{\Pn}, \quad g \mapsto g.x. \]
  When $x$ corresponds to a Borel-fixed ideal, it is fixed by $B$ and
  so by \cite[Sec. 12.1]{Hu} we get a morphism
  \[ G/B \lpil \Hilb^{\Pn}. \]
  But $G/B$ is a projective variety, \cite[Section 21.3]{Hu}, hence complete and
  so the image, the
  orbit of $x$, is a closed subvariety of $\Hilb^{\Pn}$,
  \cite[Section 21.1]{Hu}.
  
  Conversely suppose an orbit $Y$ of the action \eqref{app:degen} is closed,
  and let $Y$ have the reduced scheme structure. We get a morphism
  $G \times Y \pil \Hilb^{\Pn}$ which factors through $Y$
  (since it is reduced)
  to give $G \times Y \pil Y$. The restriction $B \times Y \pil Y$
  has a fixed point by the Borel fixed point theorem \cite[Section 21.2]{Hu}.
  This fix point corresponds to an ideal $I$ such that $g.I = I$ for
  every $g \in B$. So $I$ is a Borel-fixed ideal, and $Y$ is its orbit
  by $G = GL(n+1)$.

  Part b). Let $I$ be a strongly stable ideals, and suppose $J = g.I$
  is also strongly stable. We show that $I$ and $J$ are equal.
  Given $g$, for each $\langle x_1, \ldots, x_i \rangle$ let $\tau(i)$ be
  minimal such that $g.\langle x_1, \ldots, x_i \rangle  \sus
  \langle x_1, \ldots, x_{\tau(i)} \rangle$.
  Then $\tau(i) \geq i$. Let $S = \{ i \, | \, \tau(i) = i \}$.
  Clearly $n \in S$. If $\tau^\prime$ is the associated function to $g^{-1}$,
  it is clear that the associated $S^\prime$ must equal $S$.

  Suppose now first $S = \{ n \}$. We show that $I$ and $J$ are both
  the ideal $x_n^d$ for some $d$. Let $m =  \prod_{p = 1}^n x_p^{i_p}$
  be a minimal strongly stable generator for $I$.
  Since $\tau(p) > p$,
  for each $p < n$ there is a $q = q(p) \leq p$ such that the $g(x_q)$ has a
  variable with index $>p$ if $p <n$ and index $n$ if $p = n$.
  Let this  index be $r = r(p)$.
  Then $\prod_{p = 1}^n x_{q(p)}^{i_p}$ is in $I$ by it being strongly stable,
  and $\ell = \prod_{p = 1}^n x_{r(p)}^{i_p}$ is in $J$, since $J$ is monomial.
  If $m$ is not a power of $x_n$, we note that $m >_{st} \ell$.
  By applying the same argument to $g^{-1}$ and $\ell$
  we get an element $m^\prime$ in $I$ with $\ell \geq m^\prime$.
  But if $m$ is not a power of $x_n$, this contradicts $m$ being minimal.
  So $m = x_n^d$ for some $d$, and so also $x_n^d$ is in $J$, and these
  ideals must be equal (they contain all monomials of degree $d$).
 
  Now consider the general case $S = \{ s_1 < s_2 < \cdots < s_r = n\}$.
  We claim that every minimal strongly stable generator of $I$ has
  the form $\prod_{u = 1}^r x_{s_u}^{j_u}$. Then $J = g.I$ must also have
  these as generators and so $J = g.I$.
  Let $m$ be a mimimal generator for $I$ and write $m = \prod_{u = 1}^r m_u$,
  where the variables in $m_u$ are $x_j$ with $s_{u-1} < j \leq s_u$.
  By the same type of argument as in the $S = \{n\}$ case, we will have
  $m^\prime = \prod_{u = 1}^r x_{s_u}^{j_u}$ in $I$ where $d_u = \deg m_u$.
  \end{proof}


  \begin{corollary} (char.$\,k = 0$)
  If $I$ is Borel-fixed and $>$ any term ordering, then
  for $g$ in an open subset of $GL(n+1)$ the initial ideal
  (the generic intial ideal) $\ini_{>}(g.I) = I$.
  \end{corollary}

  \begin{proof}
    The initial ideal $\ini_>(g.I)$ is the limit at $t = 0$ of a family of
  ideals parametrized by $\Spec k[t]$, \cite[Section 3.2]{He-Hi}
  and whose general member is
  a coordinate change of $g.I$ and so of $I$. Thus $\ini_>(g.I)$ is in
  the closure of the orbit of $g.I$ and so is {\it in} the orbit of $g.I$
  or equivalently of $I$. But since the generic initial ideal
  $\ini_>(g.I)$ is Borel-fixed \cite[Chap.15]{Ei},
  it is then equal to $I$.
\end{proof}

\begin{remark}
  Several people have to their surprise observed the above corollary.
  It is stated and shown
  for $g$ in an open subset of $GL(n+1)$ in \cite[Prop. 4.2.6(b)]{He-Hi}, and
  they attribute it to A.Conca. Theorem \ref{thm:degen} has been known
  by M.Stillman
  since the late 1980's, who learned it from D.Bayer. He also informed that
  the theorem and its corollary  essentially follow from Borel's fixed point
  theorem, as shown above.

  Galligo's theorem \cite{Gal} that any ideal degenerates to a Borel ideal, and
  the theorem above are inspiration for approaches to the classification
  of Hilbert scheme components using Borel ideals, \cite{Ber-Lel-Rog},
  \cite{Cio-Rog, Fl-Ro, Ka-Le}, and recently \cite{Ram-bor, Sta-hil}.
\end{remark}
  

\bibliographystyle{amsplain}
\bibliography{biblio}
\end{document}